%% file: main.tex
\documentclass{article}
\usepackage[margin=1in]{geometry}
\usepackage{setspace}

\usepackage{amsmath,amssymb,amsthm,mathtools}
\usepackage{graphicx}
\usepackage{subcaption}
\usepackage{tikz}
\usepackage{xcolor}
\usepackage{natbib}

\theoremstyle{plain}
\newtheorem{theorem}{Theorem}[section]
\newtheorem{proposition}[theorem]{Proposition}
\newtheorem{corollary}[theorem]{Corollary}
\newtheorem{lemma}[theorem]{Lemma}

\newtheorem{assumption}[theorem]{Assumption}
\theoremstyle{definition}
\newtheorem{definition}[theorem]{Definition}
\theoremstyle{remark}

\newcommand{\abs}[1]{\left\lvert#1\right\rvert}

\title{Minimax-Optimal Static Estimators under Interference}
\author{Arthur Delarue$^1$, Ilan Lobel$^2$}
\date{{\footnotesize$^1$UVA Darden School of Business, \texttt{delaruea@darden.virginia.edu}\\$^2$NYU Stern School of Business, \texttt{ilobel@stern.nyu.edu}}\\\vspace{10pt}
October 31, 2026}
\begin{document}

\maketitle

\begin{abstract}
    We study the problem of estimating the global treatment effect in randomized experiments under interference, where a unit's outcome may depend on the treatment assignments of other units. We adopt a minimax framework, allowing an adversary to select an outcome function from a class of functions that satisfy bounded outcomes and bounded interference. We derive minimax-optimal static estimators that are robust to a wide range of interference structures, parameterized by a scalar $\delta$ controlling the maximum influence of one unit's treatment on another's outcome. Under large interference, we observe two regimes. When the treatment and control groups are equally sized, the minimax-optimal static estimator is the zero estimator. When one of the two groups is larger, the minimax-optimal static estimator is a ``half-zero'' estimator, in which the outcomes of the smaller group are ignored, and the estimator only considers the outcomes of the larger group.
\end{abstract}

\section{Introduction}

Digital companies often run randomized experiments to evaluate the impact of new features on user behavior and key business metrics. The simplest experiment is an A/B test, where users are randomly assigned to either a treatment group that receives the new feature or a control group that does not. The goal is to estimate the global treatment effect, defined as the difference in outcomes between the state of the world where all users receive the treatment and the state of the world where all users receive the control. A core assumption for valid inference in this design is that a user's outcome depends only on their own treatment assignment and not on the treatment assignments of other users. This assumption is often violated in practice, as users may interact with each other and influence each other's outcomes, in a phenomenon known as interference.

Typical approaches to mitigate interference follow one of two big ideas. The first idea is to design the experiment to reduce interactions between treated and control users, with techniques like cluster randomization or switchback designs. The second idea is to modify the estimator to account for an assumed interference structure. Examples include the exposure model, which assumes that a user's outcome depends on their own treatment assignment and the number of treated neighbors, and the linear-in-means model, which assumes that a user's outcome depends linearly on their own treatment assignment and the average treatment assignment of their neighbors.

In this paper, we study the second category of approaches through a robust optimization lens. Instead of assuming a specific model of interference, we allow an adversary to select an outcome function from a class of functions that satisfy two natural constraints: bounded outcomes and bounded interference. We then seek an estimator that minimizes the worst-case mean-squared error (in expectation over an experiment in which each of $n$ units is independently assigned to treatment with probability $p$) over this class of outcome functions. This minimax approach allows us to derive estimators that are robust to a wide range of interference structures, all parametrized by a single scalar $\delta$ that controls the maximum amount by which the outcome of one unit can be affected by the treatment assignment of another. We restrict our attention to so-called ``static'' estimators, in which each observed outcome is weighted only according to its own treatment or control status. This class of estimators includes the classical Horvitz-Thompson estimator, which ignores interference.

Our analysis of this minimax model yields several key insights. We first consider the no-interference setting, in which $\delta=0$, and show that the classical Horvitz-Thompson estimator is suboptimal for finite $n$. The minimax-optimal estimator adds a novel shrinkage factor, which reduces the worst-case relative risk from $O(n^{-1/2})$ to $O(n^{-1})$. Reassuringly, the shrinkage factor vanishes as $n\to\infty$, so the Horvitz-Thompson estimator is asymptotically minimax-optimal.

Our second, more significant, contribution establishes a threshold value of the interference parameter $\delta$ beyond which the minimax-optimal estimator is the zero estimator, which ignores all outcomes and gives up on estimation altogether. When the experiment is symmetric (i.e., $p=1/2$), this threshold is the point at which the collective effect of a unit's neighbors' treatments is on par with the unit's own treatment, so that the experimenter cannot distinguish between the two. This threshold provides a clear recommendation for practitioners: if the magnitude of interference is too large, a static estimator is simply the wrong data-driven tool to support decision-making.

Our analysis of asymmetric experiments ($p \neq 1/2$) sheds further light on this finding. In this case, we find that there are two thresholds. When $\delta$ exceeds the first one, the minimax-optimal static estimator only considers the larger group (treatment if $p > 1/2$, control if $p < 1/2$) and sets the estimated value of the smaller group to zero. When $\delta$ exceeds the second threshold, the minimax-optimal static estimator is again the zero estimator. For asymptotically large $n$, this means that the experimenter can \emph{decrease} worst-case risk by running a maximally imbalanced experiment (e.g., $p=1$).  Our model thus suggests a robust estimation strategy is to estimate outcomes under global treatment directly (which can be done without any interference), then use assumptions to bound outcomes under global control.

\section{Related Work}

The question of interference has received increased attention in the operations literature because of the ubiquity of A/B testing in digital platforms. One line of work focuses on improving the experiment design to reduce interference, including randomizing treatment assignments over clusters of users \citep{ugander2013graph,holtz2025reducing}, over time \citep{bojinov2023design}, or both \citep{ni2023design,jia2023clustered}.

Another line of work focuses on improving the estimator to account for interference, usually by assuming a particular operational model and the interference structure it implies. Building on the general model of \cite{aronow2017estimating}, recent approaches develop estimators that account for queueing \citep{li2026experimenting}, capacity constraints \citep{boutilier2024operational}, matching marketplaces \citep{bright2025reducing} and inventory stockouts \citep{chen2025b}. \cite{farias2022markovian} relate estimation under interference to off-policy learning in dynamic programming and develop an estimator that leverages optimal Q functions to reduce bias. Similarly, \cite{peng2025differences} develop a ``difference-in-neighbors'' estimator to reduce estimation error using an underlying network structure.

Other approaches seek higher-level managerial insights about the impact of interference on downstream business decisions. \cite{johari2024does} show that interference can sometimes increase statistical power without affecting false positive rates. \cite{han2024population} study panel experiments in which treatment is sequentially randomized; they find that observing outcomes over time can help reduce interference, but these gains can be mitigated or erased by the presence of carryover effects.

Our work considers the problem of estimation under interference through a robust optimization lens \citep[see, e.g.,][for thorough reviews of robust optimization principles]{ben2009robust,bertsimas2022robust}. We assume that potential outcomes lie in a particular uncertainty set, and we seek a robust estimator that minimizes the worst-case mean-squared error over this set. Robust optimization formulations have been used recently to find minimax-optimal randomization designs, for switchback experiments \citep{bojinov2023design} and clustered experiments \citep{brennan2022cluster,candogan2024correlated,viviano2023causal}. \cite{thiyageswaran2026optimal} seek to design a minimax-optimal clustered design to balance interference bias reduction with covariate balance. \cite{kandiros2026design} are the first to study minimax-optimal estimators: their work develops upper and lower bounds for the minimax risk obtained from selecting both design and estimator, without explicitly characterizing the estimator. Though their goal is to find a minimax design, our work is close to \cite{viviano2023causal}: in particular, we adopt a very similar class of outcome functions for the inner maximization problem. \cite{sudijono2026sharp} study minimax-optimal estimators in a different setting without interference, similarly finding a minimax-optimal shrinkage estimator, though with different assumptions and a different shrinkage factor.

We note that robustness has a slightly different meaning in the statistics literature, where it is often used to describe estimators that are insensitive to outliers or model misspecification. We refer to robustness in the optimization sense, i.e., as the ability of an estimator to perform well under a range of possible outcome functions. \cite{blanchet2025distributionally} provide a detailed discussion of these two notions of robustness.

\section{A Minimax Model of Interference}\label{sec:model}

\subsection{Model Setup}

We consider an undirected $d$-regular graph $G=(\mathcal{V}, \mathcal{E})$, denoting the number of nodes $|\mathcal{V}|$ by $n$ and writing $\mathbf A$ for its adjacency matrix, so that $\mathbf A_{k\ell} = \mathbb I[(k,\ell)\in\mathcal E]$. Note that because each node has $d$ neighbors, we assume $nd$ is even, and we additionally assume $n > 1$. We assume that $G$ is \emph{vertex-transitive}: for every pair of nodes $i,j\in\mathcal V$, there is a relabeling of the nodes that maps $i$ to $j$. Formally, for every pair $i,j\in \mathcal{V}$, there exists a bijection $\pi:\mathcal V\to\mathcal V$ satisfying $(k,\ell)\in\mathcal E \iff (\pi(k),\pi(\ell))\in\mathcal E$ for all $k,\ell\in\mathcal V$, such that $\pi(i)=j$. We call such an edge-preserving bijection an \emph{automorphism} of $G$, and write $\mathrm{Aut}(G)$ for the set of them. A vertex-transitive graph ``looks the same'' from the perspective of every node: examples include the complete graph and any cycle.

The binary variable $z_i$ denotes whether node $i$ (also called unit $i$) is assigned to treatment $(z_i=1)$ or control $(z_i=0)$. An outcome function $y=(y_i)_{i=1}^n$ specifies an outcome $y_i(\mathbf z)$ for every unit $i$ and assignment $\mathbf z\in\{0,1\}^n$. We write $\mathbf y(\mathbf z)=(y_i(\mathbf z))_{i=1}^n$ for the resulting outcome vector.
Our objective is to estimate the global treatment effect, defined as
\[
GTE(y) = \sum_{i=1}^n y_i(\mathbf{1}) - \sum_{i=1}^n y_i(\mathbf{0}).
\]

We cannot simultaneously observe the global control and global treatment states. Instead, we can sample a treatment assignment for each unit, according to $Z_i\sim \text{Bernoulli}(p)$, then use the assignment and observed outcomes to estimate the global treatment effect. Our goal is to minimize the worst-case mean-squared error over a particular class of outcome functions. Let $\mathcal{Y}(\mathcal{E}, \delta)$ denote the class of outcome functions satisfying the following two assumptions.

\begin{assumption}[Bounded outcomes]\label{assum:bounded-outcomes}
For every unit $i$ and treatment assignment $\mathbf z\in\{0,1\}^n$, $y_i(\mathbf z)\in[-1,1]$.
\end{assumption}

\begin{assumption}[Bounded interference]\label{assum:bounded-interference}
For every treatment assignment $\mathbf z\in\{0,1\}^n$ and every pair of distinct units $i,j$,
\[
    \left|y_i(\mathbf{z}^{(z_j=1)}) - y_i(\mathbf{z}^{(z_j=0)})\right| \le \delta \mathbf A_{ij},
\]
where $\mathbf{z}^{(z_j=k)}$ denotes the treatment assignment vector obtained by setting the $j$-th entry of $\mathbf{z}$ to value $k$.
\end{assumption}

Assumption~\ref{assum:bounded-outcomes} is self-explanatory. Assumption~\ref{assum:bounded-interference} ensures that changing the treatment assignment of unit $j\neq i$ affects the outcome of unit $i$ by at most $\delta$ if $i$ and $j$ share an edge in the underlying graph, and does not affect the outcome of unit $i$ if $i$ and $j$ do not share an edge.

%For tractability, we will also sometimes consider a restricted class of \emph{exchangeable} outcome functions. Let $N(i)$ designate unit $i$'s neighbors, and $s_i(\mathbf{z})$ the number of its neighbors that are treated, i.e., $s_i(\mathbf{z}) = \sum_{j\in N(i)} z_j$. We say an outcome function is exchangeable if there exists a function $\bar{y}: \{0,1\}\times\{0,\ldots,d\} \to [-1,1]$, satisfying $\left|\bar{y}(z_i, k+1) - \bar{y}(z_i, k)\right| \le \delta$ for all $z_i \in \{0,1\},\ k \in \{0,\ldots,d-1\}$, such that $y_i(\mathbf z) = \bar{y}(z_i, s_i(\mathbf z))$ for each unit $i$. The set of outcome functions in $\mathcal{Y}(\mathcal{E}, \delta)$ that satisfy this additional restriction is denoted by $\mathcal{Y}^{\text{exch}}(\mathcal E,\delta)$.

We consider the class of static linear estimators, $a=(\mathbf a^0,\mathbf a^1)$ consisting of two vectors of fixed weights, with $a_i^1$ applied to unit $i$ under treatment and $a_i^0$ under control. For an assignment $\mathbf z$, write $\mathbf a(\mathbf z)=(a_i^1z_i+a_i^0(1-z_i))_{i=1}^n$; the resulting estimate is $\mathbf a(\mathbf z)^{\mathsf T}\mathbf y(\mathbf z)$. Let $\mathcal A$ denote the class of these estimators. Given a treatment assignment probability $p$, the worst-case risk (or mean-squared error) of estimator $a$ is
\[
F_p(a)\coloneqq\max_{y\in\mathcal Y(\mathcal E,\delta)}
\mathbb E_{\mathbf Z\sim\operatorname{Bernoulli}(p)^n}
\left[\left(\mathbf a(\mathbf Z)^{\mathsf T}\mathbf y(\mathbf Z)-GTE(y)\right)^2\right].
\]
When $p$ is fixed, we write $F(a)$ for $F_p(a)$. For weights that are constant across units, we also write $F(a^0,a^1)$ for $F((a^0\mathbf 1,a^1\mathbf 1))$.

Our goal is then to solve the following minimax problem:
\begin{equation}\label{eq:minimax}
\min_{a\in\mathcal{A}} \max_{y\in \mathcal{Y}(\mathcal{E}, \delta)}\mathbb{E}_\mathbf{Z}\left[\left(\mathbf a(\mathbf Z)^{\mathsf T}\mathbf y(\mathbf Z) - GTE(y)\right)^2\right].%,
\end{equation}
%where $\mathcal{A}$ denotes either $\mathcal{A}_{\text{static}}$ or $\mathcal{A}_{\text{adaptive}}$.

Equation~\ref{eq:minimax} models interference as an adversarial process, in which the adversary's strength is controlled by the single parameter $\delta$. When $\delta=0$, there is no interference at all; as $\delta$ increases, the adversary can select outcome functions that are increasingly sensitive to the treatment assignments of neighboring units.

Our focus on static estimators will allow us to significantly reduce the dimensionality of the outer minimization problem and thus derive strong properties on the optimal solution and its associated worst-case risk. While the class of static estimators is restricted, it does include the classical Horvitz-Thompson estimator, which we recall below.

\begin{definition}[Horvitz--Thompson estimator]\label{def:ht}
For $p\in(0,1)$, the Horvitz--Thompson estimator $a_{HT}\in\mathcal A$ has the same weights for every unit:
\[
a_{HT}^1=\frac1p, \qquad a_{HT}^0=-\frac1{1-p}.
\]
We write $F^{HT}(p)\coloneqq F_p(a_{HT})$ for its worst-case risk.
\end{definition}

Before further analyzing this minimax problem, we observe that assuming bounded outcomes is necessary for~\eqref{eq:minimax} to be well-defined.

\begin{proposition}
    \label{prop:unbounded-outcomes}
If outcomes are unbounded, the minimax mean-squared error in~\eqref{eq:minimax} is infinite for any $\delta$.
\end{proposition}

The proof of Proposition~\ref{prop:unbounded-outcomes} is straightforward: if outcomes are unbounded, then the adversary can select an outcome function that affects only one unit, with arbitrarily large treatment effect. Because that single unit is either treated or not, the estimator cannot observe both potential outcomes, and the mean-squared error is also arbitrarily large.

\subsection{Estimator Symmetry}

We first show that the minimax problem~\eqref{eq:minimax} can be reduced to a lower-dimensional problem. Because the underlying graph is vertex-transitive, the optimal estimator can be chosen to have constant weights across units.

\begin{theorem}[Symmetric optimal estimator]\label{thm:sym-static-general-p}
Let $G$ be an undirected graph and $p\in(0,1)$. If $\delta=0$ or $G$ is vertex-transitive, there is an optimal solution $(\mathbf a^0,\mathbf a^1)$ to problem~\eqref{eq:minimax} with $a_i^0=a^0$ and $a_i^1=a^1$ for every $i$, for some scalars $a^0,a^1\in\mathbb R$.
\end{theorem}

Theorem~\ref{thm:sym-static-general-p} allows us to consider only two weights, $a^0$ and $a^1$, rather than $2n$ weights. For general $p$, the treatment and control weights may differ because the two groups appear in different proportions. When $p=1/2$, a further reduction is possible.

\begin{corollary}\label{cor:sym-static-vt}
Let $p=1/2$ and suppose $\delta=0$ or $G$ is vertex-transitive. There is an optimal solution to problem~\eqref{eq:minimax} of the form $a_i(\mathbf z)=a(2z_i-1)$ for every $i$, for some scalar $a\in\mathbb R$.
\end{corollary}

Indeed, when treatment and control are assigned in equal proportions, the optimal estimator can be chosen to be symmetric, with $a^1=a=-a^0$ for some scalar $a$. For such symmetric static estimators, we use the further scalar shorthand $F(a)\coloneqq F(-a,a)$ for $a\in\mathbb R$, in a slight abuse of notation.

In the symmetric case, we can also show that the treatment unit weight is always non-negative, which further restricts problem~\eqref{eq:minimax} to the non-negative real line when $p=1/2$.

\begin{proposition}[Non-negativity of the treatment weight]\label{prop:a-nonneg}
Let $G$ be vertex-transitive and $p=1/2$. Then $\min_{a\in\mathbb R} F(a) = \min_{a\ge0} F(a)$.
\end{proposition}

Proposition~\ref{prop:a-nonneg} allows us to restrict the outer minimization in~\eqref{eq:minimax} to $a\ge0$ when $p=1/2$.

\section{Minimax-Optimal Estimators Without Interference}

We begin by studying our minimax framework in the absence of interference to understand how our analysis aligns with classical results in statistics. We first study the worst-case risk of the standard Horvitz--Thompson estimator. When $\delta=0$, let $y_i^1,y_i^0\in[-1,1]$ be unit $i$'s outcomes under treatment and control. The Horvitz--Thompson estimator is unbiased, so its mean squared error is its variance. Let $\mathcal R(a,y)$ denote the risk of estimator $a$ under outcome function $y$. We can write
\[
\mathcal R(a_{HT},y)
=p(1-p)\sum_{i=1}^n\left(\frac{y_i^1}{p}+\frac{y_i^0}{1-p}\right)^2
\le\frac{n}{p(1-p)}.
\]
Equality holds when $y_i^1=y_i^0=1$ for every $i$, hence $F^{HT}(p)=\frac{n}{p(1-p)}$. A natural question is whether the Horvitz--Thompson estimator is minimax-optimal in our framework. We show that it is only asymptotically optimal, and that it is outperformed by a novel shrinkage estimator for finite $n$.

We first consider a symmetric experiment ($p=1/2$). In this case, we can solve problem~\eqref{eq:minimax} in closed form, and obtain the following minimax-optimal estimator.

\begin{theorem}[Optimal estimator, $p=1/2$]\label{thm:no-interference-exact}
Let $\delta=0$ and $p=1/2$. The symmetric estimator
\[
a^1 = \frac{2(n-\sqrt n)}{n - 1} = -a^0
\]
is minimax-optimal, with worst-case risk
$F^\star = \frac{4n^2}{(\sqrt n + 1)^2}$.
\end{theorem}

The minimax-optimal estimator therefore shrinks the Horvitz--Thompson weights by a factor of $\frac{n-\sqrt n}{n-1}=\frac{\sqrt n}{\sqrt n+1}$. This shrinkage factor tends to one as $n$ grows large, and the Horvitz-Thompson estimator is thus asymptotically optimal. More precisely, we can quantify its relative excess risk as
\[
\frac{F^{HT}(1/2) - F^*}{F^*} = \frac{2}{\sqrt{n}} + O\left(\frac{1}{n}\right),
\]
which tends to zero as $n$ tends to infinity.

For general $p\in(0,1)$, obtaining the minimax-optimal estimator in closed form requires solving a quartic equation, and no simple solution is available. Shrinking the Horvitz--Thompson weights by the same factor nevertheless yields a near-optimal estimator.

\begin{proposition}[Near-optimal estimator, general $p$]\label{prop:no-interference-general-p}
Let $\delta = 0$, fix $p\in(0,1)$, and let $F^\star(p) \coloneqq \min_{a\in\mathcal A}F_p(a)$. The static estimator
\[
\hat a^1 = \frac1p \cdot \frac{\sqrt n}{\sqrt n+1}, \qquad \hat a^0 = -\frac1{1-p}\cdot\frac{\sqrt n}{\sqrt n+1}
\]
satisfies, as $n\to\infty$,
\[
\frac{F(\hat a) - F^\star(p)}{F^\star(p)} \;\le\; \frac{(1-2p)^2}{4p(1-p)}\cdot\frac1n + O(n^{-3/2}).
\]
\end{proposition}

Because no closed form is available for the optimal risk, we compute a Lagrangian bound to establish the result in Proposition~\ref{prop:no-interference-general-p}. We can then use this result and the closed form risks for the Horvitz--Thompson and shrinkage Horvitz--Thompson estimators to obtain the following relative excess risk result.

\begin{corollary}\label{prop:ht-excess}
Let $\delta=0$ and fix $p\in(0,1)$. As $n\to\infty$,
\[
\frac{F^{HT}(p)-F^\star(p)}{F^\star(p)}=\frac{2}{\sqrt n}+O(n^{-1}).
\]
\end{corollary}

The Horvitz-Thompson estimator's worst-case relative risk exceeds the optimum by a factor $O(n^{-1/2})$, compared to a factor $O(n^{-1})$ for the shrinkage Horvitz--Thompson estimator. The difference is meaningful for small to moderate $n$ but vanishes as $n$ tends to infinity.

Thus, our minimax framework based on robust optimization largely mirrors classical results in statistics, which establish the optimality of the Horvitz--Thompson estimator under mild distributional assumptions on outcomes. This convergence of results suggests considering estimator optimality from an adversarial rather than a stochastic perspective is reasonable when the environment exhibits no interference between units. To our knowledge, the shrinkage factor $\sqrt{n}/(\sqrt{n}+1)$ is novel and may be of independent interest.

\section{Large Interference and Experimentation Failure}

Our minimax framework largely reproduces known results in the absence of interference, at least when $n$ is large. We next consider the opposite regime, where interference is quite large relative to outcomes. Our main result is that when interference is too large, the minimax-optimal estimator is 0. We first study the special case where $p=1/2$.

\subsection{Symmetric Experiments}

We establish a sufficient condition for optimality of the zero estimator in the following theorem.

\begin{theorem}\label{thm:give-up-d-regular}
Let $G$ be vertex-transitive with common degree $d$, and $p=1/2$. If $\delta \ge \delta_d \coloneqq 2/d$, then the zero estimator, with $a_i(\mathbf z) = 0$ for every $i$ and every $\mathbf z$, is minimax-optimal, with worst-case risk $4n^2$. Moreover, it is the unique minimax-optimal symmetric static estimator: $F(a) > F(0) = 4n^2$ for every $a\neq0$.
\end{theorem}

Theorem~\ref{thm:give-up-d-regular} establishes a sufficient condition for which no estimator improves upon the worst-case risk of the zero estimator. The threshold $\delta_d=2/d$ is intuitive: under Assumption~\ref{assum:bounded-outcomes} (bounded outcomes), the maximum treatment effect for a single unit is 2. Thus, when each neighbor's treatment assignment can affect a unit's outcome by $2/d$, treating all of a unit's neighbors together can have a comparable effect as treating the unit itself. Under such conditions, the adversary can render any estimator worse than the zero estimator.

The proof of Theorem~\ref{thm:give-up-d-regular} is also instructive. We achieve the minimax bound by examining a particular outcome function which gives every estimator higher risk than the zero estimator. This outcome function is defined as $y_i(\mathbf{z})=1-\delta_d s_i(\mathbf{z})$, where $s_i(\mathbf{z})$ designates the number of treated neighbors of unit $i$ under assignment $\mathbf{z}$. To understand why this outcome function makes estimation particularly difficult, we consider a simple example. On the complete graph $K_5$, we have $d=4$ and $\delta_d=1/2$. If $K=\sum_i Z_i$ units are treated, then every treated unit has $s_i=K-1$ treated neighbors and every control unit has $s_i=K$ treated neighbors. Figure~\ref{fig:giveup-k5-example} shows the resulting outcomes and estimates for each possible number of treated units.

\begin{figure}[htbp]
\centering
\begin{minipage}[c]{0.20\textwidth}
\centering
\begin{tikzpicture}
  \coordinate (v1) at (90:1.20cm);
  \coordinate (v2) at (18:1.20cm);
  \coordinate (v3) at (-54:1.20cm);
  \coordinate (v4) at (-126:1.20cm);
  \coordinate (v5) at (162:1.20cm);
  \foreach \i/\j in {1/2,1/3,1/4,1/5,2/3,2/4,2/5,3/4,3/5,4/5}
    \draw[gray!65, line width=0.5pt] (v\i)--(v\j);
  \foreach \i in {1,...,5}
    \node[draw,circle,fill=white,inner sep=0pt,minimum size=5.5mm] at (v\i) {$\i$};
\end{tikzpicture}
\end{minipage}\hfill
\begin{minipage}[c]{0.79\textwidth}
\centering
\small
\setlength{\tabcolsep}{2pt}
\renewcommand{\arraystretch}{1.25}
\begin{tabular}{@{}c c c c c c c r@{\hspace{2pt}}c@{\hspace{2pt}}l r@{}}
\hline
& & \multicolumn{2}{c}{Treated} & \multicolumn{2}{c}{Control} & & & & & \\
\cline{3-4}\cline{5-6}
$k$ & $\Pr(K=k)$ & $s_i$ & $y_i^1$ & $s_i$ & $y_i^0$ & $y_i^1-y_i^0$ & & \makebox[0pt][c]{$\hat\tau_{HT}$} & & $\hat\tau_0$ \\
\hline
0 & $1/32$ & --- & --- & $0$ & $1$ & --- & $2(0-5\cdot1)$ & $=$ & $-10$ & $0$ \\
1 & $5/32$ & $0$ & $1$ & $1$ & $1/2$ & $1/2$ & $2(1\cdot1-4\cdot\tfrac12)$ & $=$ & $-2$ & $0$ \\
2 & $10/32$ & $1$ & $1/2$ & $2$ & $0$ & $1/2$ & $2(2\cdot\tfrac12-3\cdot0)$ & $=$ & $2$ & $0$ \\
3 & $10/32$ & $2$ & $0$ & $3$ & $-1/2$ & $1/2$ & $2(3\cdot0-2\cdot(-\tfrac12))$ & $=$ & $2$ & $0$ \\
4 & $5/32$ & $3$ & $-1/2$ & $4$ & $-1$ & $1/2$ & $2(4\cdot(-\tfrac12)-1\cdot(-1))$ & $=$ & $-2$ & $0$ \\
5 & $1/32$ & $4$ & $-1$ & --- & --- & --- & $2(5\cdot(-1)-0)$ & $=$ & $-10$ & $0$ \\
\hline
\end{tabular}
\end{minipage}
\caption{The complete graph $K_5$ at $p=1/2$ and $\delta_d=1/2$. Each row groups assignments with $k$ treated units. For every mixed assignment, the observed between-arm difference $y^1-y^0=1/2$ has the opposite sign from the global treatment effect, $GTE=-10$.}
\label{fig:giveup-k5-example}
\end{figure}

In this example, moving from the global control state to the global treatment state changes all units' outcomes from $1$ to $-1$, which means the global treatment effect is $-2n=-10$ (maximally negative under bounded outcomes). However, as long as the numbers of treated and control units are both positive, the outcome for treated units is always larger than the outcome for control units. As a result, we observe that the Horvitz--Thompson estimator (and in fact, any static estimator) grossly overestimates the treatment effect. The potential interference bias is so strong that defaulting to the zero estimator reduces worst-case risk.

The $\delta_d$ threshold thus provides an actionable guideline for practitioners running symmetric experiments. If they anticipate that the magnitude of interference can exceed this threshold, a randomized experiment is likely not the best tool to estimate the global treatment effect, and it may be worthwhile to investigate alternative decision support methodologies.

\subsection{General Experiments}

We now extend the analysis to general assignment probabilities $p\in(0,1)$. For any $r\in(0,1)$, let $B_r\sim\operatorname{Bin}(d,r)$, and define $\Psi_r(\delta)\coloneqq\mathbb E\left[\max\{-1,1-\delta B_r\}\right]$. We then define the following two thresholds:
\[
\underline\delta(p)\coloneqq\frac1{d\max(p,1-p)},\qquad \bar\delta(p)\coloneqq\inf\{\delta\ge0:\Psi_{\min(p,1-p)}(\delta)\le0\}\ge \frac1{d\min(p,1-p)}.
\]

The last inequality is proved in the appendix. Using these two thresholds, we can establish the following generalization of Theorem~\ref{thm:give-up-d-regular}.

\begin{theorem}\label{thm:give-up-general-p}
Let $G$ be vertex-transitive with $n$ vertices and degree $d\ge1$, and let $p\in(0,1)$.
\begin{itemize}
    \item For each fixed $p\ne1/2$, there are constants $C_p,N_p>0$ depending only on $p$ such that, for $n\ge N_p$ and $\delta\ge(1+C_p/n)\underline\delta(p)$, there exists a minimax-optimal symmetric static estimator that does not use the minority group: $a^0=0$ if $p>1/2$, and $a^1=0$ if $p<1/2$.
    \item The zero estimator is minimax-optimal if and only if $\delta\ge\bar\delta(p)$.
\end{itemize}
\end{theorem}

\begin{figure}[h!]
    \centering
    \vspace{-10pt}
    \includegraphics{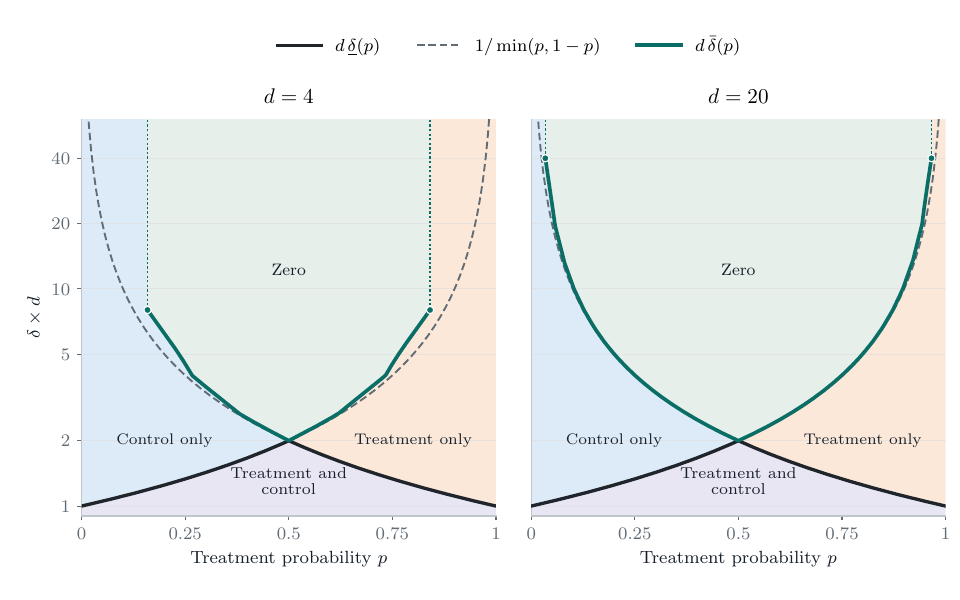}
    \vspace{-10pt}
    \caption{Illustration of Theorem~\ref{thm:give-up-general-p} for $d=4$ and $d=20$. The curve $d\,\underline\delta(p)=1/\max(p,1-p)$ is an asymptotic sufficient bound for dropping the minority group; the dashed curve $1/\min(p,1-p)$ is a lower bound on the exact zero-estimator threshold $d\,\bar\delta(p)$. For fixed $p\ne1/2$, an optimal estimator using only the majority group exists above the black curve once $n$ is sufficiently large. The treatment-and-control shading below $d\underline\delta(p)$ is illustrative because the exact switch is unknown. The zero estimator is optimal for sufficiently large $\delta$ if and only if $\max(p,1-p)^d\le1/2$.}
    \label{fig:giveup-thresholds}
\end{figure}

Theorem~\ref{thm:give-up-general-p} shows that the single give-up threshold of the symmetric experiment splits into two distinct thresholds. For sufficiently large $n$ and $\delta\ge(1+C_p/n)\underline\delta(p)$, a minimax-optimal static estimator can be chosen to set the weight of the minority group to zero. When $\delta\ge\bar\delta(p)$, giving up on estimation altogether is minimax-optimal.

We establish in the appendix that the bound $\bar\delta(p)\ge1/[d\min(p,1-p)]$ is strict unless $p=1/2$, in which case both thresholds are equal ($\bar{\delta}(p) = \underline\delta(p)$) and we recover Theorem~\ref{thm:give-up-d-regular}. We can also show that $\bar\delta(p)$ is finite if and only if $\max(p,1-p)^d\le1/2$. Thus, when $p$ is large enough (or small enough), the zero estimator is never minimax-optimal.

Figure~\ref{fig:giveup-thresholds} illustrates the thresholds in Theorem~\ref{thm:give-up-general-p} for $d=4$ and $d=20$. For small $\delta$, the optimal estimator has nonzero weights for both treatment and control units. As $\delta$ increases, the minimax-optimal estimator first sets the weight of the smaller group (e.g., control if $p > 1/2$) to zero. Eventually, it becomes optimal to set both weights to zero if $\max(p,1-p)^d\le1/2$; otherwise, the zero estimator never becomes optimal, even as $\delta$ grows to infinity.

From a practice perspective, both thresholds are relevant. Once $\delta$ sufficiently exceeds $\underline\delta(p)$ (and $n$ is large enough), the minimax-optimal estimator completely ignores the outcomes of the smaller group. If the experiment has a small treatment and large control group, this means ignoring the treatment outcomes altogether. Considering the optimal risk as a function of $p$ and $\delta$ allows us to refine this takeaway. We first define the minimax-optimal risk as
\[
V_{G,p}(\delta)\coloneqq\min_{a\in\mathcal A}F_p(a).
\]
We then establish the following bounds on this optimal risk.

\begin{theorem}[Finite-sample bounds on optimal risk]\label{thm:risk-bounds}
Let $G$ be vertex-transitive with $n$ vertices and degree $d\ge1$, and let $p\in(0,1)$. Define
\begin{equation}\label{eq:risk-curve}
\Phi_{d,p}(\delta)\coloneqq
\big[\min\{dp\delta,1\}+\min\{d(1-p)\delta,1\}\big]^2.
\end{equation}
For every $\delta\ge0$,
\begin{equation}\label{eq:risk-bounds}
\max\!\left\{0,\Phi_{d,p}(\delta)-\frac{p^2+(1-p)^2}{dp(1-p)}\right\}
\le\frac{V_{G,p}(\delta)}{n^2}
\le\min\!\left\{4,\Phi_{d,p}(\delta)+\frac{4}{np(1-p)}\right\}.
\end{equation}
\end{theorem}

For fixed $p$, the upper bound converges to $\min(4, \Phi_{d,p}(\delta))$ as $n$ tends to infinity, and the lower bound converges to $\max(0, \Phi_{d,p}(\delta))$ as $d$ tends to infinity. In the limit as both $n$ and $d$ tend to infinity, the optimal risk (normalized by $n^2$) converges to $\Phi_{d,p}(\delta)$ (appropriately clipped to remain inside $[0,4]$).

\begin{figure}[t]
    \centering
    \begin{subfigure}[t]{0.49\textwidth}
        \centering
        \includegraphics[width=\linewidth]{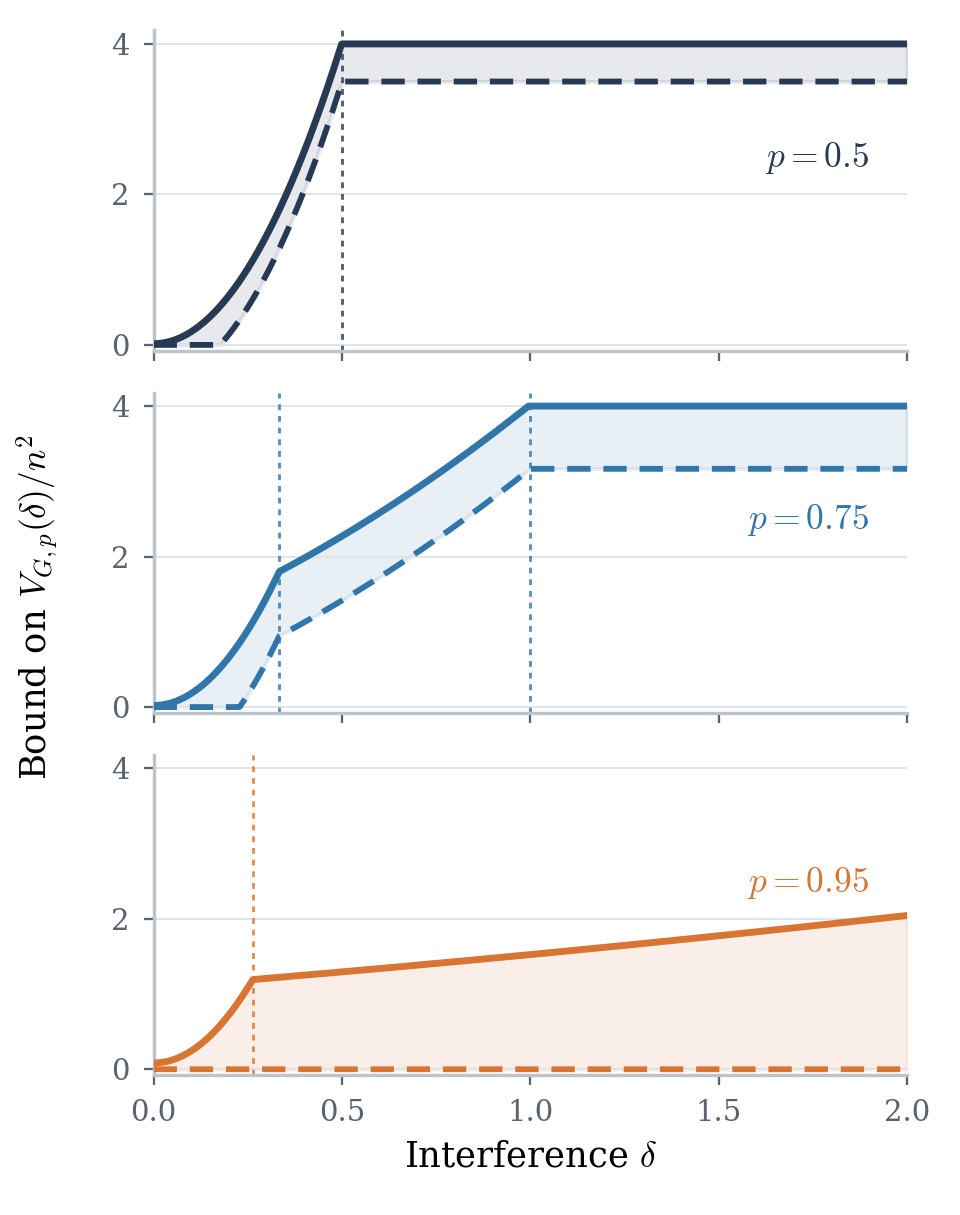}
        \caption{$d=4$}
        \label{fig:risk-bounds-d4}
    \end{subfigure}\hfill
    \begin{subfigure}[t]{0.49\textwidth}
        \centering
        \includegraphics[width=\linewidth]{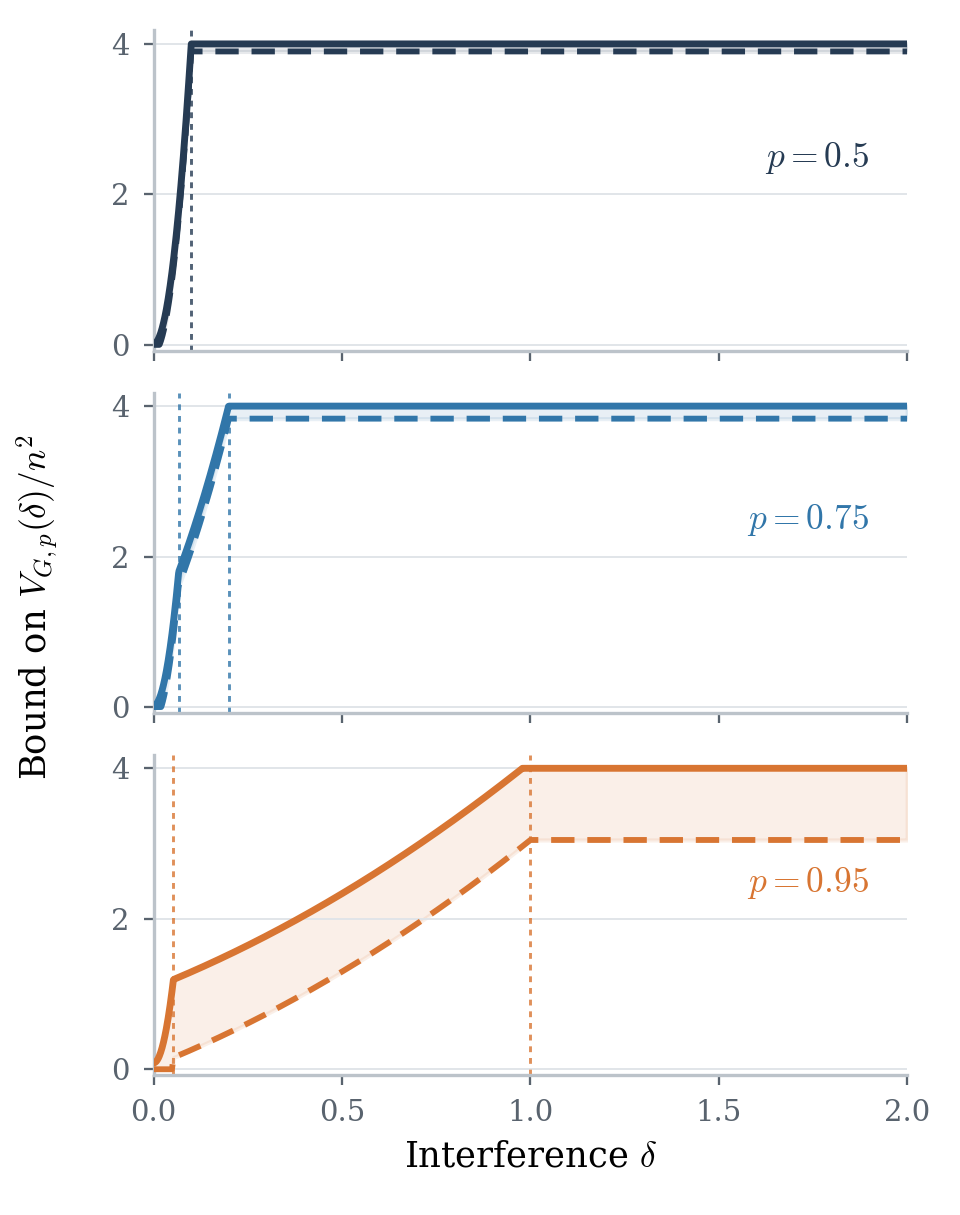}
        \caption{$d=20$}
        \label{fig:risk-bounds-d20}
    \end{subfigure}
    \caption{Bounds on normalized minimax risk $V_{G,p}(\delta)/n^2$ from Theorem~\ref{thm:risk-bounds} at $n=1000$, for $p\in\{0.5,0.75,0.95\}$ and $\delta\in[0,2]$. Solid curves show the upper bound and dashed curves the lower bound in~\eqref{eq:risk-bounds}. Vertical dotted lines mark the breakpoints $1/(dp)$ and $1/[d(1-p)]$ of $\Phi_{d,p}$.}
    \label{fig:risk-bounds}
\end{figure}

Figure~\ref{fig:risk-bounds} shows the lower and upper bounds in Theorem~\ref{thm:risk-bounds} as a function of $\delta$ for $n=1000$, $d\in\{4,20\}$, and three assignment probabilities. We observe regime changes determined by $\Phi_{d,p}(\delta)$. Under balanced assignment, $\Phi_{d,p}(\delta)=\min\{d\delta,2\}^2$: it grows quadratically until $\delta=2/d$, then becomes constant. For $p>1/2$, the curve has three regions:
\begin{equation}\label{eq:risk-regimes}
\Phi_{d,p}(\delta)
=\begin{cases}
d^2\delta^2,&0\le\delta\le1/(dp),\\
(1+d(1-p)\delta)^2,&1/(dp)\le\delta\le1/[d(1-p)],\\
4,&\delta\ge1/[d(1-p)].
\end{cases}
\end{equation}
The quadratic slope drops significantly at $\delta=1/(dp)$, and growth stops at $\delta=1/(d(1-p))$. We observe that the lower bound is not particularly tight for small degree $d$ and large treatment probability $p$. In addition, the second threshold $1/(d(1-p))$ is not always meaningful: because outcomes are bounded in $[-1,1]$, the set of outcome functions is the same for all $\delta\ge 2$. When $d=4$ and $p=0.95$, we have $1/(d(1-p))=5 > 2$: as a result, the bottom plot in Figure~\ref{fig:risk-bounds-d4} only exhibits one threshold when $\delta\in[0,2]$.

\begin{figure}[t]
    \centering
    \includegraphics{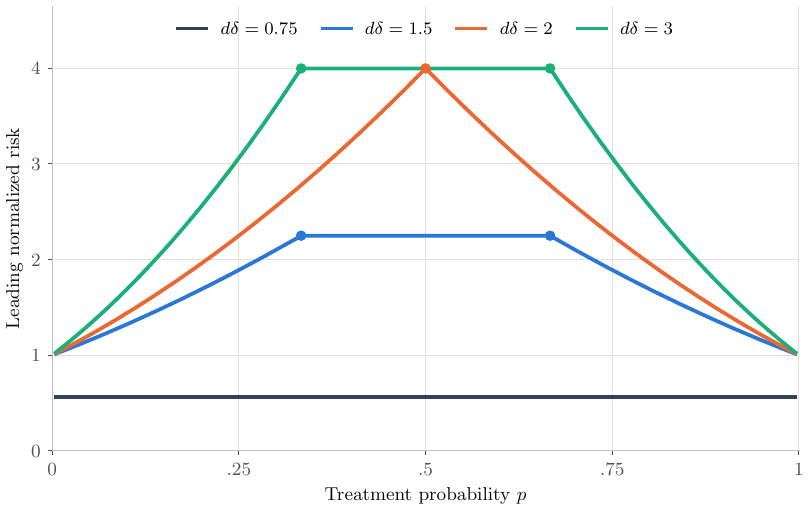}
    \caption{Optimal-risk large-$n$ upper bound $\Phi_{d,p}(\delta)$ as a function of $p$, for fixed $d\delta$. Moving away from $p=1/2$ weakly reduces the upper bound, with a strict reduction once the assignment leaves the central plateau.}
    \label{fig:risk-by-assignment}
\end{figure}

One surprising observation from Figure~\ref{fig:risk-bounds} is that for constant $\delta$, the worst-case risk seems to be non-increasing (and sometimes, sharply decreasing) as $p$ moves away from $1/2$. To further probe this mechanism, we plot $\Phi_{d,p}(\delta)$ as a function of $p$ for fixed $d\delta$ in Figure~\ref{fig:risk-by-assignment}. The curve $\Phi_{d,p}(\delta)$ is non-increasing as $p$ moves away from $1/2$. If $d\delta\le1$, it equals $d^2\delta^2$ for every $p$. Once $d\delta>1$, sufficient asymmetry strictly improves the upper bound on risk.

Figure~\ref{fig:risk-by-assignment} has a rather clear takeaway: if interference is strong enough (i.e., $\delta \ge 1/d$), the experimenter is better off running a maximally asymmetric experiment ($p\to 0$ or $p\to 1$). Indeed, this corresponds to directly estimating only global treatment (or only global control) and relying on the generic outcome bounds to gauge the value of the other state. The benefit of this approach is that estimating global treatment (or global control) directly is interference-free by construction. This view usefully complements the takeaway from Theorem~\ref{thm:give-up-general-p}: once $\delta \ge 1/d$, we will obtain a lower worst-case risk by estimating only one half of the global treatment effect (control or treatment) rather than estimating both. In this setting, interference makes estimation via an A/B test and a static estimator worse than simply rolling out the treatment, measuring outcomes, and comparing to prior knowledge (e.g., bounds) about the global control state.

\bibliographystyle{plainnat}
\bibliography{references}

\input{appendix}

\end{document}

%% file: appendix.tex
% !TeX root = ./model.tex

\newpage

\appendix

\section{Bounded outcomes and interference}\label{app:boundedness}

\begin{proof}[Proof of Proposition~\ref{prop:unbounded-outcomes}]
    Assume outcomes are unbounded, and consider the outcome function satisfying $y_1(\mathbf{z}) = Mz_1$ for some large $M>0$ and $y_i(\mathbf{z}) = 0$ for all $i \neq 1$. This outcome function is feasible for any $\delta \geq 0$, and has $GTE(y) = M$. For any static estimator $a$, we have $\mathbf a(\mathbf z)^{\mathsf T}\mathbf y(\mathbf z) = a_1^1 z_1 M$.
    Therefore, we can write the risk as
    \[
        \mathcal{R}(a, y) = \mathbb{E}_{\mathbf{Z}} \left[ (a_1^1 Z_1 M - M)^2 \right] = M^2 \mathbb{E}_{\mathbf{Z}} \left[ (a_1^1 Z_1 - 1)^2 \right] = M^2\left[p (a_1^1 - 1)^2 + (1-p) \right] \ge (1-p)M^2,
    \]
    which diverges to infinity as $M \to \infty$. Hence, the minimax risk is unbounded.
\end{proof}

\section{Symmetry reductions}\label{app:symmetry-reductions}

As in Section~\ref{sec:model}, $y$ denotes an outcome function and $a$ an estimator; bold symbols denote vectors or matrices. Thus $\mathbf y(\mathbf z)$ and $\mathbf a(\mathbf z)$ are the corresponding outcome and weight vectors under assignment $\mathbf z$, and
\[
\mathbf a(\mathbf z)^{\mathsf T}\mathbf y(\mathbf z) = \sum_{i=1}^n a_i(\mathbf z)\,y_i(\mathbf z), \qquad GTE(y)  = \sum_{i=1}^n\big(y_i(\mathbf 1)-y_i(\mathbf 0)\big) = \mathbf 1^{\mathsf T}\big(\mathbf y(\mathbf 1)-\mathbf y(\mathbf 0)\big).
\]
A static estimator $a$ is determined by its two weight vectors $\mathbf a^0,\mathbf a^1$, through $\mathbf a(\mathbf z) = \mathbf a^0 + \operatorname{diag}(\mathbf a^1-\mathbf a^0)\,\mathbf z$. Finally, we write $\mathcal R(a,y) \coloneqq \mathbb E_{\mathbf Z}\big[(\mathbf a(\mathbf Z)^{\mathsf T}\mathbf y(\mathbf Z)-GTE(y))^2\big]$, so that $F(a) = \max_{y\in\mathcal Y(\mathcal E,\delta)}\mathcal R(a,y)$. Recall that $F$ is convex, as a maximum of convex functions of $a$.

For a permutation $\pi$ of $[n]$, let $\mathbf P_\pi$ be the associated permutation matrix such that $\mathbf P_\pi \mathbf e_k = \mathbf e_{\pi(k)}$, or equivalently $(\mathbf P_\pi \mathbf v)_i = v_{\pi^{-1}(i)}$ for any vector $\mathbf v\in \mathbb{R}^n$. Recall the following properties of permutation matrices:
\[
\mathbf P_\pi\mathbf P_\tau = \mathbf P_{\pi\tau}, \qquad \mathbf P_\pi^{-1} = \mathbf P_\pi^{\mathsf T}, \qquad \mathbf P_\pi\mathbf 1 = \mathbf 1, \qquad \mathbf P_\pi^{\mathsf T}\mathbf e_j = \mathbf e_{\pi^{-1}(j)}.
\]
We write $\mathbf P$ for $\mathbf P_\pi$ when $\pi$ is fixed. In this notation $\pi$ is an automorphism of $G$ exactly when $\mathbf P^{\mathsf T}\mathbf A\mathbf P = \mathbf A$, equivalently $\mathbf P\mathbf A\mathbf P^{\mathsf T} = \mathbf A$.

We define the relabeled estimator $\pi\cdot a$ and the relabeled outcome function $\pi\cdot y$ by
\[
\pi\cdot a \coloneqq \big(\mathbf P\mathbf a^0,\ \mathbf P\mathbf a^1\big), \qquad\qquad (\pi \cdot \mathbf{y})(\mathbf{z}) \coloneqq \mathbf P\,\mathbf y(\mathbf P^{\mathsf T}\mathbf z).
\]
Using $\operatorname{diag}(\mathbf P\mathbf v) = \mathbf P\operatorname{diag}(\mathbf v)\mathbf P^{\mathsf T}$ and $\mathbf P^{\mathsf T}\mathbf P = \mathbf I$, we can write, for every $\mathbf z$,
\begin{equation}\label{eq:static-relabel}
(\pi\cdot\mathbf a)(\mathbf P\mathbf z) = \mathbf{P}\mathbf a^0 + \operatorname{diag}(\mathbf P(\mathbf a^1-\mathbf a^0))\,\mathbf P\mathbf z = \mathbf P\big(\mathbf a^0 + \operatorname{diag}(\mathbf a^1-\mathbf a^0) \mathbf{P}^{\mathsf T}\mathbf{P}\mathbf z\big) = \mathbf P\,\mathbf a(\mathbf z).
\end{equation}

\begin{lemma}[Relabeling preserves the risk]\label{lem:relabel-identity}
Let $a\in\mathcal A$, let $y$ be any outcome function, and let $\pi$ be any permutation of $[n]$. Then $\mathcal R(\pi\cdot a,\pi\cdot y) = \mathcal R(a,y)$ for every $p\in(0,1)$. If in addition $\pi\in\mathrm{Aut}(G)$, or if $\delta=0$, then $y\in\mathcal Y(\mathcal E,\delta)\iff\pi\cdot y\in\mathcal Y(\mathcal E,\delta)$, and hence $F(\pi\cdot a)=F(a)$.
\end{lemma}

\begin{proof}
\emph{Risk.} For every $\mathbf z$, by~\eqref{eq:static-relabel} and the definition of $\pi\cdot y$,
\[
(\pi\cdot\mathbf a)(\mathbf P\mathbf z)^{\mathsf T}(\pi\cdot\mathbf y)(\mathbf P\mathbf z) = \big(\mathbf P\,\mathbf a(\mathbf z)\big)^{\mathsf T}\mathbf P\,\mathbf y(\mathbf z) = \mathbf a(\mathbf z)^{\mathsf T}\mathbf y(\mathbf z),
\]
the last step by $\mathbf P^{\mathsf T}\mathbf P=\mathbf I$. The $GTE$ is unchanged as well: using $\mathbf P^{\mathsf T}\mathbf 1=\mathbf 1$, $\mathbf P^{\mathsf T}\mathbf 0=\mathbf 0$ and $\mathbf 1^{\mathsf T}\mathbf P=\mathbf 1^{\mathsf T}$,
\[
GTE(\pi\cdot y) = \mathbf 1^{\mathsf T}\mathbf P\big(\mathbf y(\mathbf 1)-\mathbf y(\mathbf 0)\big) = GTE(y).
\]
Because the $Z_i$ are i.i.d., permuting the coordinates of $\mathbf Z$ leaves its distribution unchanged: $\mathbf P\mathbf Z \stackrel{d}{=}\mathbf Z$. We may therefore compute the risk of $(\pi\cdot a,\pi\cdot y)$ at $\mathbf P\mathbf Z$ in place of $\mathbf Z$, and then apply the two identities above:
\[
\mathcal R(\pi\cdot a,\pi\cdot y) = \mathbb E_{\mathbf Z}\Big[\big((\pi\cdot\mathbf a)(\mathbf P\mathbf Z)^{\mathsf T}(\pi\cdot\mathbf y)(\mathbf P\mathbf Z)-GTE(\pi\cdot y)\big)^2\Big] = \mathbb E_{\mathbf Z}\Big[\big(\mathbf a(\mathbf Z)^{\mathsf T}\mathbf y(\mathbf Z)-GTE(y)\big)^2\Big] = \mathcal R(a,y).
\]

\emph{Feasibility.} Boundedness is immediate. For bounded interference, let $\mathbf w=\mathbf P^{\mathsf T}\mathbf z$; then $\mathbf P^{\mathsf T}\mathbf e_j = \mathbf e_{\pi^{-1}(j)}$ gives $\mathbf P^{\mathsf T}\big(\mathbf z^{(z_j=1)}\big) = \mathbf w^{(w_{\pi^{-1}(j)}=1)}$, so
\[
\abs{(\pi\cdot y)_i\big(\mathbf z^{(z_j=1)}\big) - (\pi\cdot y)_i\big(\mathbf z^{(z_j=0)}\big)} = \abs{y_{\pi^{-1}(i)}\big(\mathbf w^{(w_{\pi^{-1}(j)}=1)}\big) - y_{\pi^{-1}(i)}\big(\mathbf w^{(w_{\pi^{-1}(j)}=0)}\big)} \le \delta\,\mathbf A_{\pi^{-1}(i)\,\pi^{-1}(j)},
\]
by feasibility of $y$. The right-hand side can be re-written as $\delta\,\mathbf A_{\pi^{-1}(i)\,\pi^{-1}(j)} = \delta\,\big(\mathbf P\mathbf A\mathbf P^{\mathsf T}\big)_{ij}$ which equals $\delta\mathbf A_{ij}$ when $\pi$ is an automorphism of $G$ or $\delta = 0$. Hence $y\mapsto\pi\cdot y$ maps $\mathcal Y(\mathcal E,\delta)$ onto itself, and maximizing the risk identity over $y$ gives $F(\pi\cdot a)=F(a)$.
\end{proof}

When $p=1/2$, we introduce $\sigma$ as the operation that exchanges the roles of treatment and control, and we define the flipped estimator $\sigma\cdot a$ and the flipped outcome function $\sigma\cdot y$ by
\[
\sigma\cdot a \coloneqq \big({-\mathbf a^1},\ {-\mathbf a^0}\big), \qquad\qquad (\sigma\cdot\mathbf y)(\mathbf z) \coloneqq -\,\mathbf y(\mathbf 1-\mathbf z).
\]
Similarly to~\eqref{eq:static-relabel}, for every $\mathbf z$, we have
\begin{equation}\label{eq:static-flip}
(\sigma\cdot\mathbf a)(\mathbf 1-\mathbf z) = -\,\mathbf a(\mathbf z).
\end{equation}

\begin{lemma}[Flipping preserves the risk]\label{lem:flip-identity}
Let $p=1/2$, let $a\in\mathcal A$, and let $y$ be any outcome function. Then $\mathcal R(\sigma\cdot a,\sigma\cdot y)=\mathcal R(a,y)$ and $y\in\mathcal Y(\mathcal E,\delta)\iff\sigma\cdot y\in\mathcal Y(\mathcal E,\delta)$, and hence $F(\sigma\cdot a)=F(a)$.
\end{lemma}

\begin{proof}
\emph{Risk.} For every $\mathbf z$, by~\eqref{eq:static-flip} and the definition of $\sigma\cdot y$, the two sign changes cancel:
\[
(\sigma\cdot\mathbf a)(\mathbf 1-\mathbf z)^{\mathsf T}(\sigma\cdot\mathbf y)(\mathbf 1-\mathbf z) = \big(-\mathbf a(\mathbf z)\big)^{\mathsf T}\big(-\mathbf y(\mathbf z)\big) = \mathbf a(\mathbf z)^{\mathsf T}\mathbf y(\mathbf z),
\]
and $GTE(\sigma\cdot y) = \mathbf 1^{\mathsf T}\big(-\mathbf y(\mathbf 0)+\mathbf y(\mathbf 1)\big) = GTE(y)$. When $p=1/2$, $\mathbf Z$ and $1 - \mathbf Z$ have the same distribution ($\mathbf 1-\mathbf Z\stackrel{d}{=}\mathbf Z$). We may therefore compute the risk of $(\sigma\cdot a,\sigma\cdot y)$ using $\mathbf 1-\mathbf Z$, and then apply the two identities above to rewrite the estimator and global treatment effect terms:
\[
\mathcal R(\sigma\cdot a,\sigma\cdot y) = \mathbb E_{\mathbf Z}\Big[\big((\sigma\cdot\mathbf a)(\mathbf 1-\mathbf Z)^{\mathsf T}(\sigma\cdot\mathbf y)(\mathbf 1-\mathbf Z)-GTE(\sigma\cdot y)\big)^2\Big] = \mathbb E_{\mathbf Z}\Big[\big(\mathbf a(\mathbf Z)^{\mathsf T}\mathbf y(\mathbf Z)-GTE(y)\big)^2\Big] = \mathcal R(a,y).
\]

\emph{Feasibility.} Boundedness is again immediate. Since $\mathbf 1-\big(\mathbf z^{(z_j=1)}\big) = (\mathbf 1-\mathbf z)^{(z_j=0)}$,
\[
(\sigma\cdot y)_i\big(\mathbf z^{(z_j=1)}\big) - (\sigma\cdot y)_i\big(\mathbf z^{(z_j=0)}\big) = y_i\big((\mathbf 1-\mathbf z)^{(z_j=1)}\big) - y_i\big((\mathbf 1-\mathbf z)^{(z_j=0)}\big),
\]
again at most $\delta\,\mathbf A_{ij}$ in absolute value by feasibility of $y$. Hence $y\mapsto\sigma\cdot y$ maps $\mathcal Y(\mathcal E,\delta)$ onto itself, and maximizing the risk identity over $y$ gives $F(\sigma\cdot a)=F(a)$.
\end{proof}

\begin{proof}[Proof of Theorem~\ref{thm:sym-static-general-p}]
Let $\Gamma$ be the group of all permutations of $[n]$ if $\delta=0$, and $\mathrm{Aut}(G)$ otherwise. In either case, $\Gamma$ acts transitively on the units. Let $a^\star$ be optimal. By Lemma~\ref{lem:relabel-identity}, $F(\pi\cdot a^\star)=F(a^\star)$ for every $\pi\in\Gamma$, so convexity makes
\[
\bar a\coloneqq\frac1{|\Gamma|}\sum_{\pi\in\Gamma}\pi\cdot a^\star
\]
optimal as well. For $k\in\{0,1\}$, its weight vector is $\bar{\mathbf a}^k=|\Gamma|^{-1}\sum_{\pi\in\Gamma}\mathbf P_\pi\mathbf a^{\star,k}$. For any $\tau\in\Gamma$, composition by $\tau$ permutes the elements of $\Gamma$, so
\[
\mathbf P_\tau\bar{\mathbf a}^k
=\frac1{|\Gamma|}\sum_{\pi\in\Gamma}\mathbf P_{\tau\pi}\mathbf a^{\star,k}
=\bar{\mathbf a}^k.
\]
Because $\Gamma$ is transitive, any two coordinates of $\bar{\mathbf a}^k$ are equal. Thus $\bar a_i(\mathbf z)=a^1z_i+a^0(1-z_i)$ for every unit $i$.
\end{proof}

\begin{proof}[Proof of Corollary~\ref{cor:sym-static-vt}]
By Theorem~\ref{thm:sym-static-general-p} there is an optimal $a^\star$ with $\mathbf a^{\star,k} = a^{\star,k}\mathbf 1$. By Lemma~\ref{lem:flip-identity} and convexity, $\bar a \coloneqq \frac12\big(a^\star+\sigma\cdot a^\star\big)$ is optimal too. Its weights are given by $\bar{\mathbf a}^1 = \frac12\big(\mathbf a^{\star,1}-\mathbf a^{\star,0}\big) = a\mathbf 1$ and $\bar{\mathbf a}^0 = \frac12\big(\mathbf a^{\star,0}-\mathbf a^{\star,1}\big) = -a\mathbf 1$, where $a \coloneqq \frac12(a^{\star,1}-a^{\star,0})$. Hence $\bar a_i(\mathbf z) = az_i - a(1-z_i) = a(2z_i-1)$.
\end{proof}

\begin{proof}[Proof of Proposition~\ref{prop:a-nonneg}]
The outcome function $\mathbf y^\dagger(\mathbf z) \coloneqq 2\mathbf z-\mathbf 1$ is feasible for every $\mathcal E$ and every $\delta\ge0$, since $y_i^\dagger$ depends only on $z_i$. It has $GTE(y^\dagger)=2n$ and, for every $\mathbf z$, $\mathbf a(\mathbf z)^{\mathsf T}\mathbf y^\dagger(\mathbf z) = a(2\mathbf z-\mathbf 1)^{\mathsf T}(2\mathbf z-\mathbf 1) = an$, so
\[
F(a) \;\ge\; \mathcal R(a,y^\dagger) \;=\; n^2(a-2)^2 \qquad \text{for every } a\in\mathbb R.
\]
Also $F(0) = \max_y GTE(y)^2 = 4n^2$, since $|y_i(\mathbf 1)-y_i(\mathbf 0)|\le2$ forces $|GTE(y)|\le 2n$, with equality at $y^\dagger$. For $a<0$ we have $(a-2)^2>4$, so $F(a) > 4n^2 = F(0)$. Hence no $a<0$ improves on $a=0$, and $\min_{a\in\mathbb R}F(a)=\min_{a\ge0}F(a)$.
\end{proof}

\section{No interference}\label{app:no-interference}

When $\delta=0$, each outcome depends only on its own unit's treatment assignment. For every $\mathcal E$, write
\[
\mathcal Y(\mathcal E,0)=\mathcal Y_0\coloneqq
\left\{y:y_i(\mathbf z)=y_i^1z_i+y_i^0(1-z_i),\quad
 y_i^0,y_i^1\in[-1,1],\quad i\in[n]\right\}.
\]

\begin{proof}[Proof of Theorem~\ref{thm:no-interference-exact}]
When $\delta=0$, every graph reduces to $n$ disconnected nodes, which is trivially vertex-transitive. By Corollary~\ref{cor:sym-static-vt}, we may assume $a_i(\mathbf z) = a(2z_i-1)$ for every $i$.

Fix such an $a$ and $y\in\mathcal Y_0$, and write $s_i \coloneqq y_i^1+y_i^0$ and $d_i \coloneqq y_i^1-y_i^0$, so that $GTE(y)=\sum_i d_i$. Since $(2z_i-1)y_i(\mathbf z) = y_i^1z_i - y_i^0(1-z_i)$,
\[
\mathbf a(\mathbf Z)^{\mathsf T}\mathbf y(\mathbf Z) = a\sum_{i=1}^n\big(y_i^1Z_i - y_i^0(1-Z_i)\big) = a\sum_{i=1}^n s_iZ_i - a\sum_{i=1}^n y_i^0 .
\]
The $Z_i$ are independent with $p=1/2$, so the estimator has variance $\frac{a^2}{4}\sum_i s_i^2$ and mean $\frac a2\sum_i s_i - a\sum_i y_i^0 = \frac a2\sum_i d_i=\frac a2 GTE(y)$. We can therefore write the estimator's risk as
\[
\mathcal R(a,y) = \underbrace{\frac{a^2}{4}\sum_{i=1}^n s_i^2}_{\text{variance}} + \underbrace{\frac{(a-2)^2}{4}\left(\sum_{i=1}^n d_i\right)^2}_{\text{squared bias}}.
\]
For fixed $a$, this is a sum of squares of affine functions of $(y_i^0,y_i^1)_{i=1}^n$, hence convex, so its maximum over the box $[-1,1]^{2n}$ occurs at a vertex: $y_i^0,y_i^1\in\{-1,1\}$ for every $i$. At such a $y$, call unit $i$ \emph{concordant} if $y_i^1=y_i^0$, in which case $s_i=\pm2$ and $d_i=0$, and \emph{discordant} if $y_i^1=-y_i^0$, in which case $s_i=0$ and $d_i=\pm2$. For $0\le m\le n$, let $V_m$ be the set of vertices with exactly $m$ concordant units, and for $y\in V_m$ let $n_+$ and $n_-$ count the discordant units with $d_i=+2$ and $d_i=-2$, so that $n_++n_-=n-m$. Every $y\in V_m$ verifies $\sum_i s_i^2=4m$ and $\sum_i d_i = 2(n_+-n_-)$, hence
\[
\left(\sum_{i=1}^n d_i\right)^2 \le 4(n-m)^2, 
\]
with equality if and only if $n_+=0$ or $n_-=0$. Such a vertex exists in $V_m$: take $(y_i^1,y_i^0)=(1,1)$ for the concordant units and $(y_i^1,y_i^0)=(1,-1)$ for the discordant ones. Therefore
\[
\max_{y\in V_m}\mathcal R(a,y) = \frac{a^2}{4}\cdot4m + \frac{(a-2)^2}{4}\cdot4(n-m)^2 = a^2m + (a-2)^2(n-m)^2 .
\]
Every vertex lies in exactly one $V_m$, so maximizing over $m$ gives
\[
F(a) = \max_{0\le m\le n}\Big[a^2m + (a-2)^2(n-m)^2\Big],
\]
and the bracket is convex in $m$, so the maximum is attained at $m=n$ or $m=0$:
\[
F(a) = \max\big\{na^2,\; n^2(a-2)^2\big\}.
\]
By Proposition~\ref{prop:a-nonneg} we know that the optimal $a$ is nonnegative. The function $na^2$ is increasing on $[0,\infty)$ and the function $n^2(a-2)^2$ is increasing on $[2,\infty)$, so the minimizer lies in $[0,2]$. On this interval, $na^2$ increases from $0$ to $4n$ while $n^2(a-2)^2$ decreases from $4n^2$ to $0$, so the two cross exactly once and the crossing point minimizes $F$. Taking square roots, the crossing solves $\sqrt n\,a = n(2-a)$, giving
\[
a^\star = \frac{2n}{n+\sqrt n} = \frac{2\sqrt n}{\sqrt n+1}.
\]
The optimal value is
\[
F^\star = n(a^\star)^2 = \frac{4n^2}{(\sqrt n+1)^2}.
\]
\end{proof}

For general $p$ the two weights $a^0,a^1$ no longer collapse to a single scalar, but the same vertex argument still reduces the adversary to two candidate outcome functions.

\begin{lemma}[Two-profile reduction, general $p$]\label{lem:no-int-profiles}
Let $\delta=0$, let $p\in(0,1)$, and let $a$ be a static estimator with $a_i^k=a^k$ for every $i$ and $k\in\{0,1\}$. Then $F(a^0,a^1) = \max\{f_1(a^0,a^1),\, f_2(a^0,a^1)\}$, where
\[
f_1 \coloneqq p(1-p)\,n(a^1-a^0)^2 + n^2\big(pa^1+(1-p)a^0\big)^2, \qquad f_2 \coloneqq p(1-p)\,n(a^1+a^0)^2 + n^2\big(pa^1-(1-p)a^0-2\big)^2.
\]
\end{lemma}

\begin{proof}
Fix $y\in\mathcal Y_0$. Since $Z_i(1-Z_i)=0$,
\[
\mathbf a(\mathbf Z)^{\mathsf T}\mathbf y(\mathbf Z) = \sum_{i=1}^n c_i Z_i + \sum_{i=1}^n a^0y_i^0, \qquad c_i \coloneqq a^1y_i^1-a^0y_i^0,
\]
so, writing $GTE=\sum_i(y_i^1-y_i^0)$ and using the bias--variance decomposition of a sum of independent Bernoulli terms,
\[
\mathcal R(a,y) = \underbrace{p(1-p)\sum_{i=1}^n c_i^2}_{\text{variance}} + \underbrace{\Big(\sum_{i=1}^n \ell_i\Big)^2}_{\text{squared bias}}, \qquad \ell_i \coloneqq y_i^1(pa^1-1) + y_i^0\big(1+a^0(1-p)\big).
\]
For fixed $a$, the risk is a sum of squares of affine functions of $(y_i^0,y_i^1)_{i=1}^n$, hence convex, so its maximum over the box $[-1,1]^{2n}$ occurs at a vertex: $y_i^0,y_i^1\in\{-1,1\}$ for every $i$. As a result, $c_i$ and $\ell_i$ each take one of four values, corresponding to the four possible outcome pairs, or \emph{profiles}, $(y_i^1,y_i^0)\in\{(1,1),(1,-1),(-1,1),(-1,-1)\}$; write $c(s,t)$ and $\ell(s,t)$ for the values taken at profile $(s,t)$. Collecting the units by profile, with counts $n_{s,t}\coloneqq|\{i : (y_i^1,y_i^0)=(s,t)\}|$,
\[
\mathcal R(a,y) = p(1-p)\sum_{s,t}n_{s,t}\,c(s,t)^2 + \Big(\sum_{s,t}n_{s,t}\,\ell(s,t)\Big)^2 .
\]
The right-hand side is a convex function of $(n_{s,t})$, so it is maximized at a vertex of the simplex, i.e., $n_{s,t}=n$, or every unit sharing a single profile $(s,t)$. Since both the estimate and $GTE$ change sign together under $y\mapsto -y$, $\mathcal R(a,-y)=\mathcal R(a,y)$ for every $y$, so only two of the four profiles need checking: $(1,1)$, which gives $f_1$, and $(1,-1)$, which gives $f_2$.
\end{proof}

At $p=1/2$ and $a^1=-a^0=a$ we recover $f_1 = na^2$ and $f_2 = n^2(a-2)^2$, i.e., the expression for $F$ obtained in the proof of Theorem~\ref{thm:no-interference-exact}.

\begin{proof}[Proof of Proposition~\ref{prop:no-interference-general-p}]
By Theorem~\ref{thm:sym-static-general-p} we may assume $a_i^k=a^k$ for every $i$, and then Lemma~\ref{lem:no-int-profiles} gives $F(a^0,a^1)=\max\{f_1(a^0,a^1),f_2(a^0,a^1)\}$. Write $\tilde p \coloneqq p(1-p)$. At $\hat a$,
\[
\hat a^1+\hat a^0 = \frac{1-2p}{\tilde p}\cdot\frac{\sqrt n}{\sqrt n+1}, \qquad p\hat a^1-(1-p)\hat a^0-2 = \frac{-2}{\sqrt n+1},
\]
so, using $(1-2p)^2+4\tilde p=1$,
\[
f_2(\hat a) = \tilde p\, n(\hat a^1+\hat a^0)^2 + n^2\big(p\hat a^1-(1-p)\hat a^0-2\big)^2 = \big[(1-2p)^2+4\tilde p\big]\cdot\frac{n^2}{\tilde p(\sqrt n+1)^2} = f_1(\hat a),
\]
so $F(\hat a) = f_1(\hat a) = \dfrac{n^2}{\tilde p(\sqrt n+1)^2}$.

To bound $F^\star(p)=\min_a\max\{f_1,f_2\}$ from below, note that $\max\{f_1(a),f_2(a)\}\ge\lambda f_1(a)+(1-\lambda)f_2(a)$ for any $\lambda\in[0,1]$, so minimizing both sides over $a$ gives,
\[
F^\star(p) \;\ge\; LB(\lambda) \coloneqq \min_{a^0,a^1}\big[\lambda f_1(a^0,a^1)+(1-\lambda)f_2(a^0,a^1)\big].
\]
Expanding $f_1$ and $f_2$ gives
\[
\begin{aligned}
\lambda f_1+(1-\lambda)f_2 ={}& (\tilde pn+p^2n^2)(a^1)^2+2(2\lambda-1)\tilde pn(n-1)\,a^1a^0+\big(\tilde pn+(1-p)^2n^2\big)(a^0)^2\\
&-4(1-\lambda)n^2\big(pa^1-(1-p)a^0\big)+4(1-\lambda)n^2.
\end{aligned}
\]
Setting both partial derivatives to zero and dividing by $2n$ gives the linear system
\[
\begin{aligned}
(\tilde p+p^2n)\,a^1+(2\lambda-1)\tilde p(n-1)\,a^0 &= 2(1-\lambda)pn,\\
(2\lambda-1)\tilde p(n-1)\,a^1+\big(\tilde p+(1-p)^2n\big)\,a^0 &= -2(1-\lambda)(1-p)n.
\end{aligned}
\]
The diagonal coefficients are positive and, using $p^2+(1-p)^2=1-2\tilde p$, the determinant is $\tilde p\,D(\lambda)$ with $D(\lambda)\coloneqq n+4\lambda(1-\lambda)\tilde p(n-1)^2>0$. The quadratic is therefore strictly convex, and its minimizer is the unique solution of the system, which by Cramer's rule is
\[
a^1_\lambda = \frac{2(1-\lambda)n\big[1+2\lambda(1-p)(n-1)\big]}{D(\lambda)}, \qquad a^0_\lambda = -\frac{2(1-\lambda)n\big[1+2\lambda p(n-1)\big]}{D(\lambda)}.
\]
Multiplying the first equation by $a^1_\lambda$ and the second by $a^0_\lambda$ and adding shows that, at the minimizer, the quadratic terms of the expansion equal $2(1-\lambda)n^2\big(pa^1_\lambda-(1-p)a^0_\lambda\big)$. Hence
\[
LB(\lambda) = 4(1-\lambda)n^2-2(1-\lambda)n^2\big(pa^1_\lambda-(1-p)a^0_\lambda\big) = \frac{4\lambda(1-\lambda)n^2\big[n-4(1-\lambda)\tilde p(n-1)\big]}{D(\lambda)},
\]
where the last step uses $pa^1_\lambda-(1-p)a^0_\lambda=\frac{2(1-\lambda)n[1+4\lambda\tilde p(n-1)]}{D(\lambda)}$ and $D(\lambda)-(1-\lambda)n\big[1+4\lambda\tilde p(n-1)\big]=\lambda\big[n-4(1-\lambda)\tilde p(n-1)\big]$.

Now take
\[
\lambda_n \coloneqq 1-\frac{1}{4\tilde p(\sqrt n+1)},
\]
which lies in $(0,1)$ for all sufficiently large $n$, namely whenever $\sqrt n>\frac{(1-2p)^2}{4\tilde p}$. Then $4(1-\lambda_n)\tilde p=\frac1{\sqrt n+1}$, so $n-4(1-\lambda_n)\tilde p(n-1)=n-\sqrt n+1$ and $4\lambda_n(1-\lambda_n)\tilde p(n-1)^2=\lambda_n(\sqrt n-1)(n-1)$, and therefore
\[
LB(\lambda_n) = \frac{\lambda_n n^2(n-\sqrt n+1)}{\tilde p(\sqrt n+1)\big[n+\lambda_n(\sqrt n-1)(n-1)\big]} > 0.
\]
Dividing $F(\hat a)=\frac{n^2}{\tilde p(\sqrt n+1)^2}$ by $LB(\lambda_n)$, and using $(\sqrt n+1)(n-\sqrt n+1)=n^{3/2}+1$ and $(\sqrt n-1)(n-1)=n^{3/2}-n-\sqrt n+1$,
\[
\frac{F(\hat a)}{LB(\lambda_n)} = \frac{n+\lambda_n(\sqrt n-1)(n-1)}{\lambda_n(n^{3/2}+1)} = 1+\frac{\sqrt n\,\big[\sqrt n-\lambda_n(\sqrt n+1)\big]}{\lambda_n(n^{3/2}+1)}.
\]
Since $\sqrt n-\lambda_n(\sqrt n+1)=\frac1{4\tilde p}-1=\frac{(1-2p)^2}{4\tilde p}$ and $F^\star(p)\ge LB(\lambda_n)$, we conclude
\[
\frac{F(\hat a)-F^\star(p)}{F^\star(p)} \;\le\; \frac{F(\hat a)}{LB(\lambda_n)}-1 \;=\; \frac{(1-2p)^2\sqrt n}{4p(1-p)\,\lambda_n(n^{3/2}+1)} \;=\; \frac{(1-2p)^2}{4p(1-p)}\cdot\frac1n+O(n^{-3/2}),
\]
where the last step uses $\lambda_n=1-O(n^{-1/2})$.
\end{proof}

\begin{proof}[Proof of Corollary~\ref{prop:ht-excess}]
The proof of Proposition~\ref{prop:no-interference-general-p} gives $F(\hat a)=n^2/[p(1-p)(\sqrt n+1)^2]$, and hence
\[
\frac{F^{HT}(p)}{F(\hat a)}=\left(1+\frac1{\sqrt n}\right)^2.
\]
Because $F^\star(p)\le F(\hat a)$, Proposition~\ref{prop:no-interference-general-p} implies $F(\hat a)/F^\star(p)=1+O(n^{-1})$ for fixed $p$. Therefore
\[
\frac{F^{HT}(p)-F^\star(p)}{F^\star(p)}
=\left(1+\frac1{\sqrt n}\right)^2\bigl(1+O(n^{-1})\bigr)-1
=\frac2{\sqrt n}+O(n^{-1}).
\]
\end{proof}

\section{Large interference}

For each unit $i$, let $N(i)$ denote its set of neighbors and write
\[
s_i(\mathbf z)\coloneqq\sum_{j\in N(i)}z_j
\]
for the number of treated neighbors under assignment $\mathbf z$.

\begin{proof}[Proof of Theorem~\ref{thm:give-up-d-regular}]
By the argument in the proof of Proposition~\ref{prop:a-nonneg} (Appendix~\ref{app:symmetry-reductions}), the zero estimator has worst-case risk $F(0)=\max_y GTE(y)^2=4n^2$.

It remains to show that $F(a)\ge4n^2$ for every $a\in\mathcal A$, with strict inequality for every nonzero symmetric static estimator. Consider the outcome function $y_i(\mathbf z)\coloneqq1-\delta_d\,s_i(\mathbf z)$. Flipping $z_j$ for $j\in N(i)$ changes $s_i(\mathbf z)$ by exactly one, and hence $y_i(\mathbf z)$ by exactly $\delta_d$; flipping $z_j$ for $j\notin N(i)$ leaves $y_i(\mathbf z)$ unchanged; and $s_i(\mathbf z)\in[0,d]$ gives $y_i(\mathbf z)\in[1-\delta_d d,\,1]=[-1,1]$. So $y\in\mathcal Y(\mathcal E,\delta)$ whenever $\delta\ge\delta_d$.

Because $G$ is $d$-regular, $s_i(\mathbf 1)=d$ and $s_i(\mathbf 0)=0$ for every $i$, so $y_i(\mathbf 1)=-1$ and $y_i(\mathbf 0)=1$, giving $GTE(y)=\mathbf 1^{\mathsf T}\big(\mathbf y(\mathbf 1)-\mathbf y(\mathbf 0)\big)=-2n$.

Now let $a\in\mathcal A$ have weight vectors $\mathbf a^0,\mathbf a^1$. The coefficient $a_i(\mathbf Z)=a_i^1Z_i+a_i^0(1-Z_i)$ depends only on $Z_i$, while $y_i(\mathbf Z)$ depends only on $(Z_j)_{j\in N(i)}$; since $i\notin N(i)$, the two are independent. Moreover, $p=1/2$ gives $\mathbb E[s_i(\mathbf Z)]=d/2$, and hence $\mathbb E[y_i(\mathbf Z)]=1-\delta_d\,d/2=0$. Therefore
\[
\mathbb E_{\mathbf Z}\big[\mathbf a(\mathbf Z)^{\mathsf T}\mathbf y(\mathbf Z)\big] = \sum_{i=1}^n\mathbb E[a_i(\mathbf Z)]\,\mathbb E[y_i(\mathbf Z)] = 0,
\]
and the bias--variance decomposition gives
\[
\mathcal R(a,y) = \mathrm{Var}\big(\mathbf a(\mathbf Z)^{\mathsf T}\mathbf y(\mathbf Z)\big) + GTE(y)^2 = \mathrm{Var}\big(\mathbf a(\mathbf Z)^{\mathsf T}\mathbf y(\mathbf Z)\big) + 4n^2 \;\ge\; 4n^2.
\]
Since $y$ is one particular feasible outcome function, $F(a)\ge\mathcal R(a,y)\ge4n^2=F(0)$ for every $a\in\mathcal A$, so the zero estimator is minimax-optimal.

Finally, for the symmetric static estimator with scalar $a\in\mathbb R$, we have $\mathbf a(\mathbf z)^{\mathsf T}\mathbf y(\mathbf z)=a\,W(\mathbf z)$, where $W(\mathbf z)\coloneqq\sum_{i=1}^n(2z_i-1)y_i(\mathbf z)$, so $\mathrm{Var}\big(\mathbf a(\mathbf Z)^{\mathsf T}\mathbf y(\mathbf Z)\big)=a^2\,\mathrm{Var}\big(W(\mathbf Z)\big)$. To compute this variance, let $T_i\coloneqq2Z_i-1$. Because, $p=1/2$, the $T_i$ are i.i.d.\  uniform on $\{-1,1\}$. Since $s_i(\mathbf Z)=\sum_{j\in N(i)}\frac{1+T_j}2=\frac d2+\frac12\sum_{j\in N(i)}T_j$ and $\delta_d=2/d$, we have $y_i(\mathbf Z)=-\frac1d\sum_{j\in N(i)}T_j$. Each edge $\{i,j\}\in\mathcal E$ appears twice in the resulting double sum, once from each endpoint, so
\[
W(\mathbf Z) = -\frac1d\sum_{i=1}^n\sum_{j\in N(i)}T_iT_j = -\frac2d\sum_{\{i,j\}\in\mathcal E}T_iT_j.
\]
Since the $T_m$ are independent with $\mathbb E[T_m]=0$ and $T_m^2=1$, any two edges $\{i,j\},\{k,\ell\}\in\mathcal E$ satisfy
\[
\mathrm{Cov}(T_iT_j,\,T_kT_\ell) = \mathbb E[T_iT_jT_kT_\ell]-\mathbb E[T_iT_j]\,\mathbb E[T_kT_\ell] = \mathbb I\big[\{i,j\}=\{k,\ell\}\big],
\]
and therefore
\[
\mathrm{Var}\big(W(\mathbf Z)\big) = \frac4{d^2}\sum_{\{i,j\}\in\mathcal E}\ \sum_{\{k,\ell\}\in\mathcal E}\mathrm{Cov}(T_iT_j,\,T_kT_\ell) = \frac4{d^2}\,|\mathcal E| = \frac4{d^2}\cdot\frac{nd}2 = \frac{2n}d.
\]
Hence $F(a)\ge\mathcal R(a,y)=\frac{2n}d\,a^2+4n^2>4n^2$ for every $a\neq0$.
\end{proof}

\subsection{General $p$}\label{app:give-up-general-p}

Without loss of generality, assume $p\ge1/2$ and write $q=1-p$. Write $\boldsymbol\alpha=(a^0,a^1)$ for the pair of uniform coefficients. For a fixed feasible outcome function $y$, write $\mathcal R(\boldsymbol\alpha,y)$ for its mean squared error and define
\[
F_\delta(\boldsymbol\alpha)\coloneqq\max_{y\in\mathcal Y(\mathcal E,\delta)}\mathcal R(\boldsymbol\alpha,y).
\]
We suppress the subscript $\delta$ when interference is fixed. Represented by their finitely many values, the feasible outcome functions form a bounded polytope defined by linear inequalities. In particular, the set of feasible outcome functions is compact, so the maximum defining $F_\delta$ is attained.

For $r\in(0,1)$, we use $B_r\sim\operatorname{Bin}(d,r)$ and $\Psi_r(\delta)=\mathbb E\max\{-1,1-\delta B_r\}$ from the theorem. The following profiles specify outcomes as functions of the number $k$ of treated neighbors, with the unit's own treatment status held fixed:
\[
v(k)=\min\{1,-1+\delta(d-k)\},\qquad
u(k)=\max\{-1,1-\delta k\},\qquad k\in\{0,\ldots,d\}.
\]
Here $v$ is the upper treatment profile and $u$ is the lower control profile. For an assignment $\mathbf z$ with $k=s_i(\mathbf z)$ treated neighbors, a treated unit with $y_i(\mathbf1)=-1$ satisfies $y_i(\mathbf z)\le v(k)$. Similarly, a control unit with $y_i(\mathbf0)=1$ satisfies $y_i(\mathbf z)\ge u(k)$. These bounds are attained by $v$ and $u$. Under independent Bernoulli assignment, $Z_i$ is independent of $(Z_j)_{j\in N(i)}$, so $s_i(\mathbf Z)$ has distribution $\operatorname{Bin}(d,p)$. Since $d-B_p$ has distribution $\operatorname{Bin}(d,q)$, the profile definitions give
\begin{equation}\label{eq:giveup-envelope-means}
\mathbb E[v(B_p)]=-\Psi_q(\delta),\qquad
\mathbb E[u(B_p)]=\Psi_p(\delta).
\end{equation}

\begin{lemma}[Properties of $\Psi_r$]\label{lem:giveup-psi}
For $d\ge1$ and $r\in(0,1)$, the following properties hold.
\begin{enumerate}
    \setlength{\itemsep}{0pt}
    \item \emph{Monotonicity.} $\Psi_r$ is continuous and nonincreasing in $\delta$, with $\Psi_r(0)=1$. For fixed $\delta$, it is nonincreasing in $r$.
    \item \emph{Lower bound on a zero.} $\Psi_r(\delta)>0$ for $\delta<1/(dr)$, and $\Psi_r(1/(dr))=0$ if and only if $r\ge1/2$.
    \item \emph{Upper bound on a zero.} If $dr\ge1$, then $\Psi_r(2/(dr))<0$.
    \item \emph{Saturation and existence of a zero.} For $\delta\ge2$, $\Psi_r(\delta)=2(1-r)^d-1$. Consequently, $\Psi_r(\delta)\le0$ for some finite $\delta$ if and only if $(1-r)^d\le1/2$.
\end{enumerate}
\end{lemma}

\begin{proof}
Continuity and monotonicity in $\delta$ hold termwise in the binomial expectation. Monotonicity in $r$ holds because $k\mapsto\max\{-1,1-\delta k\}$ is nonincreasing and $B_r$ is stochastically increasing in $r$. Writing $x_+=\max(x,0)$,
\[
\Psi_r(\delta)=1-\delta dr+\mathbb E(\delta B_r-2)_+.
\]
The right-hand side is positive for $\delta < 1/(dr)$. It is zero at $1/(dr)$ if and only if $B_r\le 2dr$ almost surely, or equivalently $r\ge1/2$, since $\Pr(B_r=d)>0$. At $\delta=2/(dr)$, the same expression gives
\[
\Psi_r\left(\frac2{dr}\right)=-1+\frac2{dr}\mathbb E(B_r-dr)_+=-1+\frac{\mathbb E|B_r-dr|}{dr},
\]
where the second equality holds because $B_r-dr$ has mean zero, so its positive and negative parts have the same expectation. By Cauchy--Schwarz, $\mathbb E|B_r-dr|\le\sqrt{\operatorname{Var}(B_r)}=\sqrt{dr(1-r)}$, and therefore
\[
\Psi_r\left(\frac2{dr}\right)\le-1+\sqrt{\frac{1-r}{dr}}<0
\]
whenever $dr\ge1>1-r$. Finally, for $\delta\ge2$ we can write
\[
\Psi_r(\delta) = \mathbb{E}\max\{-1,1-\delta B_r\} = 1 \cdot \Pr(B_r = 0) + (-1) \cdot \Pr(B_r > 0) = 2\Pr(B_r = 0) - 1 = 2(1-r)^d - 1,
\]
which is constant and nonpositive exactly when $(1-r)^d\le1/2$. Since $\Psi_r$ is nonincreasing, this constant is its smallest value.
\end{proof}

\begin{figure}[t]
    \centering
    \includegraphics{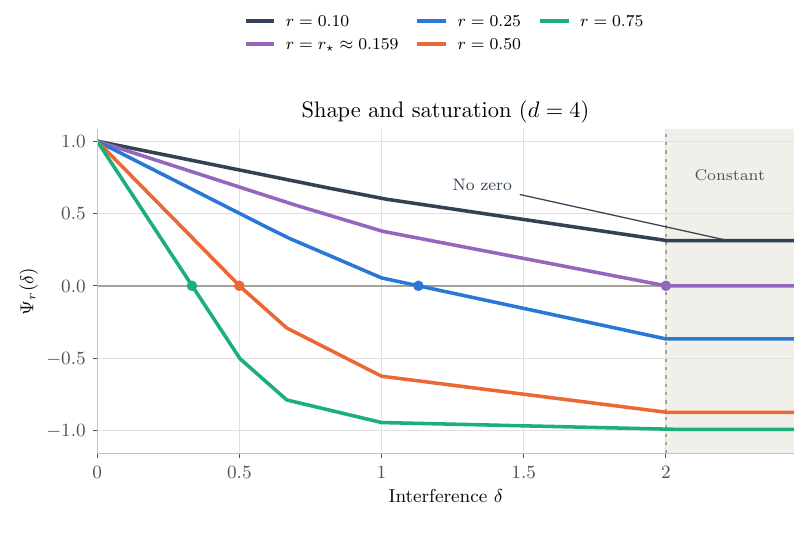}
    \caption{Monotonicity and saturation of $\Psi_r$ in Lemma~\ref{lem:giveup-psi}, illustrated at $d=4$. Each curve starts at one, decreases continuously, and becomes constant at $2(1-r)^d-1$ for $\delta\ge2$ (shaded). Increasing $r$ lowers the curve. The critical value $r_\star=1-2^{-1/d}$ separates positive terminal values, with no zero, from nonpositive terminal values; its curve first reaches zero at $\delta=2$. Dots mark first zeros.}
    \label{fig:psi-properties}
\end{figure}

Figure~\ref{fig:psi-properties} illustrates the monotonicity and saturation of $\Psi_r$. Monotonicity, positivity before $1/(dr)$, and saturation identify and compare the complete-give-up thresholds. The sign at $2/(dr)$ bounds $d\delta$ in Lemma~\ref{lem:giveup-control-profiles}.

The improvement and control-exclusion arguments below compare the worst-case risk at a candidate optimum with its value at nearby coefficient pairs. The next lemma expresses its directional derivatives using the outcome functions that attain the worst-case risk.

\begin{lemma}[Directional derivatives of the worst-case risk]\label{lem:giveup-danskin}
Fix $\delta\ge0$. For a uniform coefficient pair $\boldsymbol\alpha=(a^0,a^1)$, let $\mathcal Y^\star(\boldsymbol\alpha)\coloneqq\arg\max_{y\in\mathcal Y(\mathcal E,\delta)}\mathcal R(\boldsymbol\alpha,y)$ be the set of worst-case outcome functions. The function $F_\delta$ is convex, and for every $\boldsymbol\alpha,\mathbf w\in\mathbb R^2$,
\begin{equation}\label{eq:giveup-directional}
\begin{aligned}
F_\delta'(\boldsymbol\alpha;\mathbf w)
&\coloneqq\lim_{t\downarrow0}\frac{F_\delta(\boldsymbol\alpha+t\mathbf w)-F_\delta(\boldsymbol\alpha)}t\\
&=\max_{y\in\mathcal Y^\star(\boldsymbol\alpha)}
\nabla_{\boldsymbol\alpha}\mathcal R(\boldsymbol\alpha,y)^{\mathsf T}\mathbf w.
\end{aligned}
\end{equation}
The maximum is attained. Moreover, $\boldsymbol\alpha$ minimizes $F_\delta$ if and only if every direction $\mathbf w$ admits some $y\in\mathcal Y^\star(\boldsymbol\alpha)$ with $\nabla_{\boldsymbol\alpha}\mathcal R(\boldsymbol\alpha,y)^{\mathsf T}\mathbf w\ge0$.
\end{lemma}

\begin{samepage}
\begin{proof}
The differentiable form of Danskin's theorem \citep[Proposition~B.25]{bertsekas1999nonlinear} states the following. Let $\Theta$ be a compact subset of a Euclidean space and let $\phi:\mathbb R^m\times\Theta\to\mathbb R$. If $\phi$ and $\nabla_{\mathbf x}\phi$ are continuous, and $\phi(\cdot,\theta)$ is convex and differentiable for each $\theta$, then $\Phi(\mathbf x)=\max_{\theta\in\Theta}\phi(\mathbf x,\theta)$ is convex and
\[
\Phi'(\mathbf x;\mathbf w)=\max_{\theta\in\Theta(\mathbf x)}
\nabla_{\mathbf x}\phi(\mathbf x,\theta)^{\mathsf T}\mathbf w,
\qquad
\Theta(\mathbf x)=\arg\max_{\theta\in\Theta}\phi(\mathbf x,\theta).
\]
Apply the theorem with $\Theta=\mathcal Y(\mathcal E,\delta)$, $\mathbf x=\boldsymbol\alpha$, and $\phi=\mathcal R$, so $\Phi=F_\delta$ and $\Theta(\mathbf x)=\mathcal Y^\star(\boldsymbol\alpha)$. The set of outcomes is a polytope, therefore compact. For fixed $y$, the estimation error is affine in $\boldsymbol\alpha$, so $\mathcal R(\cdot,y)$ is a convex quadratic; its coefficients are polynomial functions of the finitely many outcome values. Thus both $\mathcal R$ and its coefficient gradient are continuous, and Danskin's theorem gives~\eqref{eq:giveup-directional}. The maximizer set $\mathcal Y^\star(\boldsymbol\alpha)$ is nonempty and compact, so the continuous gradient attains the displayed maximum.

Finally, convexity gives
\[
F_\delta(\boldsymbol\alpha+\mathbf w)-F_\delta(\boldsymbol\alpha)
\ge F_\delta'(\boldsymbol\alpha;\mathbf w).
\]
Hence $\boldsymbol\alpha$ minimizes $F_\delta$ if and only if every directional derivative is nonnegative. Equation~\eqref{eq:giveup-directional} and attainment of its maximum give the stated optimality condition.
\end{proof}
\end{samepage}

\begin{samepage}
\begin{lemma}[Optimality of zero]\label{lem:giveup-lower}
If $\Psi_q(\delta)\le0$, the zero estimator is minimax-optimal among static estimators, with risk $4n^2$. If $\Psi_q(\delta)>0$, then $F(0,t)<4n^2$ for all sufficiently small $t>0$, so the zero estimator is not minimax-optimal.
\end{lemma}
\end{samepage}

\begin{proof}
\emph{Optimality when $\Psi_q(\delta)\le0$.} Since $p\ge q$, Lemma~\ref{lem:giveup-psi} gives $\Psi_p(\delta)\le\Psi_q(\delta)\le0$. Fix a static estimator $a$ with weight vectors $\mathbf a^0,\mathbf a^1$. For each unit $i$, define the treatment and control profiles
\[
\tilde v_i(k)=\begin{cases}v(k),&a_i^1>0,\\-1,&a_i^1\le0,\end{cases}
\qquad
\tilde u_i(k)=\begin{cases}u(k),&a_i^0<0,\\1,&a_i^0\ge0.\end{cases}
\]
The outcome function
\[
y_i(\mathbf z)=z_i\tilde v_i(s_i(\mathbf z))+(1-z_i)\tilde u_i(s_i(\mathbf z))
\]
is feasible and has $GTE=-2n$. Since $i\notin N(i)$, independent treatment assignment makes $Z_i$ independent of $s_i(\mathbf Z)=\sum_{j\in N(i)}Z_j$, and hence of both $\tilde v_i(s_i(\mathbf Z))$ and $\tilde u_i(s_i(\mathbf Z))$. Also $s_i(\mathbf Z)$ has the same distribution as $B_p\sim\operatorname{Bin}(d,p)$. Thus
\begin{align*}
\mathbb E\big[Z_i\tilde v_i(s_i(\mathbf Z))\big]
&=\mathbb E[Z_i]\,\mathbb E\big[\tilde v_i(s_i(\mathbf Z))\big]
=p\,\mathbb E[\tilde v_i(B_p)],\\
\mathbb E\big[(1-Z_i)\tilde u_i(s_i(\mathbf Z))\big]
&=\mathbb E[1-Z_i]\,\mathbb E\big[\tilde u_i(s_i(\mathbf Z))\big]
=q\,\mathbb E[\tilde u_i(B_p)].
\end{align*}
Substituting the chosen outcomes into the static estimator and summing over units gives
\begin{align}
\mathbb E\big[\mathbf a(\mathbf Z)^{\mathsf T}\mathbf y(\mathbf Z)\big]
&=\sum_{i=1}^n\left\{
 a_i^1\mathbb E\big[Z_i\tilde v_i(s_i(\mathbf Z))\big]
+a_i^0\mathbb E\big[(1-Z_i)\tilde u_i(s_i(\mathbf Z))\big]\right\}\nonumber\\
&=\sum_{i=1}^n\left\{p\,a_i^1\mathbb E[\tilde v_i(B_p)]
+q\,a_i^0\mathbb E[\tilde u_i(B_p)]\right\}.\label{eq:giveup-twoterms-intermediate}
\end{align}
We can then show each term in~\eqref{eq:giveup-twoterms-intermediate} is nonnegative: if $a_i^1 > 0$, $\tilde{v}_i=v$ and $\mathbb E[\tilde v_i(B_p)]=\mathbb E[v(B_p)]=-\Psi_q(\delta)\ge0$ by~\eqref{eq:giveup-envelope-means}; if $a_i^1\le 0$, $\tilde v_i=-1$ and $\mathbb E[\tilde v_i(B_p)]=-1 \le 0$. Similarly, if $a_i^0<0$, $\tilde u_i=u$ and $\mathbb E[\tilde u_i(B_p)]=\mathbb E[u(B_p)]=\Psi_p(\delta)\le0$; if $a_i^0\ge0$, $\tilde u_i=1$ and $\mathbb E[\tilde u_i(B_p)]=1\ge0$. Consequently,
\begin{align*}
F(a)
&\ge\mathbb E\left[\big(\mathbf a(\mathbf Z)^{\mathsf T}\mathbf y(\mathbf Z)+2n\big)^2\right]\\
&=4n^2+4n\mathbb E\big[\mathbf a(\mathbf Z)^{\mathsf T}\mathbf y(\mathbf Z)\big]
+\mathbb E\big[(\mathbf a(\mathbf Z)^{\mathsf T}\mathbf y(\mathbf Z))^2\big]
\ge4n^2.
\end{align*}
The zero estimator attains this bound, so it is minimax-optimal.

\emph{Improvement when $\Psi_q(\delta)>0$.} Start with the zero estimator, whose coefficients are $(a^0,a^1)=(0,0)$ and whose worst-case risk is $4n^2$. We keep $a^0=0$ and increase $a^1$ from zero to a small positive number $t$.

An outcome function attains risk $4n^2$ at $(0,0)$ if and only if $|GTE|=2n$. Assume without loss of generality that $GTE=-2n$. Since outcomes lie in $[-1,1]$, this equality requires $y_i(\mathbf1)=-1$ and $y_i(\mathbf0)=1$ for every $i$. Treatment outcomes are then bounded above by $v(s_i(\mathbf Z))$, and~\eqref{eq:giveup-envelope-means} gives
\[
\left.\frac{d}{dt}\mathcal R((0,t),y)\right|_{t=0}
=4n\mathbb E\left[\sum_iZ_i y_i(\mathbf Z)\right] \le 4n\mathbb E\left[\sum_iZ_i v(s_i(\mathbf Z))\right]
\le-4n^2p\Psi_q(\delta)<0.
\]
Thus every outcome function that is worst-case for the zero estimator has decreasing risk when we initially increase the treatment coefficient. Reversing all outcome signs preserves the risk at every coefficient pair, so the same bound holds for worst-case outcome functions with $GTE=2n$. Equation~\eqref{eq:giveup-directional} with direction $\mathbf w=(0,1)$ therefore gives
\[
F'\big((0,0);(0,1)\big)\le-4n^2p\Psi_q(\delta)<0,
\]
so, by Lemma~\ref{lem:giveup-danskin}, $F(0,t)<F(0,0)=4n^2$ for all sufficiently small $t>0$.
\end{proof}

The remaining results prepare the control-exclusion argument. We first bound the variance of functions of the treatment assignments, then study the best estimator using only treatment outcomes, and finally construct two control outcome functions with opposite risk derivatives. We write $M_q\coloneqq1+2/q$ throughout these results.

\begin{lemma}[Variance under independent treatment assignment]\label{lem:giveup-variance}
Let $f:\{0,1\}^n\to\mathbb R$, and suppose that changing one unit's assignment changes $f$ by at most $L\ge0$. Under independent assignments $Z_i\sim\operatorname{Bernoulli}(p)$, with $q=1-p$,
\begin{equation}\label{eq:minority-bounded-differences}
\operatorname{Var}(f(\mathbf Z))\le npqL^2.
\end{equation}
\end{lemma}

\begin{proof}
The Efron--Stein inequality \citep[Theorem~3.1]{boucheron2013concentration} states that, for independent random variables $X_1,\ldots,X_n$ and a square-integrable function $W=g(X_1,\ldots,X_n)$,
\[
\operatorname{Var}(W)\le\sum_{j=1}^n\mathbb E\big[\operatorname{Var}(W\mid\mathbf X_{-j})\big],
\]
where $\mathbf X_{-j}$ denotes all coordinates except $X_j$. Apply the inequality with $X_j=Z_j$ and $W=f(\mathbf Z)$. Conditional on $\mathbf Z_{-j}$, independence leaves $Z_j$ Bernoulli with parameter $p$, so
\[
\operatorname{Var}(f(\mathbf Z)\mid\mathbf Z_{-j})
=pq\big(f(\mathbf Z^{(Z_j=1)})-f(\mathbf Z^{(Z_j=0)})\big)^2\le pqL^2.
\]
Summing over $j$ proves~\eqref{eq:minority-bounded-differences}.
\end{proof}

\begin{samepage}
\begin{lemma}[Control profiles with opposite risk derivatives]\label{lem:giveup-control-profiles}
Suppose $p>1/2$, $0\le\delta\le2$, $\Psi_q(\delta)>0$, and
\[
n>\frac{8qM_q^2}{p},\qquad c\coloneqq\frac{1+12M_q/n}{p}<2,\qquad \delta\ge c/d.
\]
Let $a\in\arg\min_{t\in\mathbb R}F(0,t)$ and consider any $y\in\mathcal Y^\star(0,a)$ with $y_i(\mathbf0)=1$ for every $i$. For every $i$ and $\mathbf z$, define
\[
\begin{aligned}
y_i^{\mathrm{lin}}(\mathbf z)&\coloneqq z_i y_i(\mathbf z)+(1-z_i)\left(1-\frac{c\,s_i(\mathbf z)}d\right),\\
y_i^{\mathrm{const}}(\mathbf z)&\coloneqq z_i y_i(\mathbf z)+(1-z_i).
\end{aligned}
\]
Then $y^{\mathrm{lin}},y^{\mathrm{const}}\in\mathcal Y^\star(0,a)$, and
\[
\left.\frac{d}{db}\mathcal R((b,a),y^{\mathrm{lin}})\right|_{b=0}<0
<\left.\frac{d}{db}\mathcal R((b,a),y^{\mathrm{const}})\right|_{b=0}.
\]
\end{lemma}
\end{samepage}

\begin{proof}
\emph{Preservation of worst-case risk.} Both constructed functions are feasible: their treatment outcomes are unchanged, their control outcomes lie in $[-1,1]$ because $0<c<2$, and changing a neighbor's assignment changes a control outcome by at most $c/d\le\delta$. Both functions also preserve $y_i(\mathbf0)=1$ and $y_i(\mathbf1)$ for every $i$, so they have the same global treatment effect as $y$. At $(0,a)$, they give the same estimate as $y$ for every assignment. Consequently,
\[
\mathcal R((0,a),y^{\mathrm{lin}})=\mathcal R((0,a),y^{\mathrm{const}})
=\mathcal R((0,a),y)=F(0,a).
\]
Thus $y^{\mathrm{lin}},y^{\mathrm{const}}\in\mathcal Y^\star(0,a)$. Write $H\coloneqq\sum_iZ_i y_i(\mathbf Z)$. For each $\widetilde y\in\{y,y^{\mathrm{lin}},y^{\mathrm{const}}\}$, the treated-outcome sum is $H$ and $GTE(\widetilde y)=GTE(y)$. Expanding the treatment derivative therefore gives the same value for all three functions:
\[
\begin{aligned}
\left.\frac{d}{dt}\mathcal R((0,t),\widetilde y)\right|_{t=a}
&=\left.\frac{d}{dt}\mathbb E\big[(tH-GTE(y))^2\big]\right|_{t=a}\\
&=2\mathbb E\big[(aH-GTE(y))H\big]\\
&=2a\mathbb E[H^2]-2GTE(y)\mathbb EH.
\end{aligned}
\]

\emph{Opposite control derivatives.} For any assignment $\mathbf z$, define
\[
T(\mathbf z)\coloneqq\sum_i z_i,\qquad
K(\mathbf z)\coloneqq\sum_i(1-z_i)\left(1-\frac{c\,s_i(\mathbf z)}d\right),
\]
and write $T\coloneqq T(\mathbf Z)$ and $K\coloneqq K(\mathbf Z)$ under randomized assignment.
Keep the treatment coefficient fixed at $a$ and vary the control coefficient $b$. The uniform estimator with weights $(b,a)$ assigns weight $a$ to treated outcomes and $b$ to control outcomes. Under the two constructed outcome functions, its values are
\[
\begin{aligned}
\sum_i[aZ_i+b(1-Z_i)]y_i^{\mathrm{lin}}(\mathbf Z)
&=a\sum_i Z_i y_i^{\mathrm{lin}}(\mathbf Z)+b\sum_i(1-Z_i)y_i^{\mathrm{lin}}(\mathbf Z)=aH+bK,\\
\sum_i[aZ_i+b(1-Z_i)]y_i^{\mathrm{const}}(\mathbf Z)
&=a\sum_i Z_i y_i^{\mathrm{const}}(\mathbf Z)+b\sum_i(1-Z_i)y_i^{\mathrm{const}}(\mathbf Z)=aH+b(n-T).
\end{aligned}
\]
Both outcome functions have global treatment effect $GTE(y)$, so their corresponding risks are
\[
\begin{aligned}
\mathcal R((b,a),y^{\mathrm{lin}})&=\mathbb E\big[(aH+bK-GTE(y))^2\big],\\
\mathcal R((b,a),y^{\mathrm{const}})&=\mathbb E\big[(aH+b(n-T)-GTE(y))^2\big].
\end{aligned}
\]
Differentiating these mean squared errors with respect to $b$ at $b=0$ gives
\begin{equation}\label{eq:giveup-control-derivatives}
\begin{aligned}
\frac12\left.\frac{d}{db}\mathcal R((b,a),y^{\mathrm{lin}})\right|_{b=0}
&=\big(a\mathbb EH-GTE(y)\big)\mathbb EK+a\operatorname{Cov}(H,K),\\
\frac12\left.\frac{d}{db}\mathcal R((b,a),y^{\mathrm{const}})\right|_{b=0}
&=\big(a\mathbb EH-GTE(y)\big)\mathbb E[n-T]+a\operatorname{Cov}(H,n-T).
\end{aligned}
\end{equation}
Here $\mathbb EK=nq(1-pc)=-12qM_q<0$ and $\mathbb E[n-T]=nq>0$. We will show that the common factor $a\mathbb EH-GTE(y)$ is at least $n$ and that the covariance terms are too small to change the resulting signs.

For the common factor, use $GTE(y)=\sum_i y_i(\mathbf1)-n$. Negating only the treatment outcomes preserves feasibility, changes the estimate to $-aH$, and changes the global treatment effect to $-GTE(y)-2n$. Since $y$ attains the worst-case risk,
\[
\mathbb E\big[(2n-aH+GTE(y))^2\big]
\le\mathbb E\big[(aH-GTE(y))^2\big].
\]
Expanding and cancelling the common squared term gives $a\mathbb EH-GTE(y)\ge n$.

To control the covariance terms, first bound $a$. Lemma~\ref{lem:giveup-lower} gives $F(0,a)<4n^2$, while the constant outcome function with all outcomes equal to one has risk $a^2\mathbb ET^2$. Hence
\[
a^2n^2p^2=a^2(\mathbb ET)^2\le a^2\mathbb E[T^2]\le F(0,a)<4n^2 \qquad \Rightarrow |a|<2/p.
\]
Next, we claim that $d\delta<2/q$: if $dq<1$, use $\delta\le2$; otherwise Lemma~\ref{lem:giveup-psi} gives $\Psi_q(2/(dq))<0<\Psi_q(\delta)$, so monotonicity implies $\delta<2/(dq)$.

To bound the change in $H$, write $H(\mathbf z)=\sum_i z_i y_i(\mathbf z)$ and fix two assignments $\mathbf z,\mathbf z'$ that differ only in the $j$-th coordinate, i.e., $z_j=0$ and $z'_j=1$. Unit $j$ contributes zero to $H(\mathbf z)$ and $y_j(\mathbf z')$ to $H(\mathbf z')$, so its contribution changes by at most one, by bounded outcomes. For every other unit $i$, the multiplier $z_i$ stays fixed. If $i\notin N(j)$, bounded interference gives $y_i(\mathbf z')=y_i(\mathbf z)$; if $i\in N(j)$, its outcome changes by at most $\delta$, and its contribution to $H$ changes only when $z_i=1$. Thus
\[
H(\mathbf z')-H(\mathbf z)
=y_j(\mathbf z')+\sum_{i\in N(j)}z_i\big(y_i(\mathbf z')-y_i(\mathbf z)\big),
\]
and the triangle inequality gives
\[
\begin{aligned}
|H(\mathbf z')-H(\mathbf z)|
&\le |y_j(\mathbf z')|+\sum_{i\in N(j)}z_i|y_i(\mathbf z')-y_i(\mathbf z)|\\
&\le 1+\delta\sum_{i\in N(j)}z_i
\le 1+d\delta<M_q,
\end{aligned}
\]
where the last two inequalities use $|N(j)|=d$ and $d\delta<2/q$, respectively. Reversing the change in $z_j$ gives the same absolute bound. Similarly, we find that $|K(\mathbf z')-K(\mathbf z)| \le 1+ c < 3$ and $|(n-T(\mathbf z'))-(n-T(\mathbf z))| \le 1$. Lemma~\ref{lem:giveup-variance} and the Cauchy--Schwarz inequality therefore give
\[
|a\operatorname{Cov}(H,K)|\le\frac2p(3npqM_q)=6nqM_q,\qquad
|a\operatorname{Cov}(H,n-T)|\le\frac2p(npqM_q)=2nqM_q.
\]
Substituting these bounds into~\eqref{eq:giveup-control-derivatives} proves the required signs:
\begin{align*}
\frac12\left.\frac{d}{db}\mathcal R((b,a),y^{\mathrm{lin}})\right|_{b=0}
&\le-12nqM_q+6nqM_q<0,\\
\frac12\left.\frac{d}{db}\mathcal R((b,a),y^{\mathrm{const}})\right|_{b=0}
&\ge n^2q-2nqM_q>0.
\end{align*}
The last inequality uses $n>4M_q$, which follows from $n>8qM_q^2/p$ because $2qM_q=2q+4>p$.
\end{proof}

\begin{proof}[Proof of Theorem~\ref{thm:give-up-general-p}]
\emph{Optimality of zero estimator.} Assume $p\ge1/2$ without loss of generality, so $q\le p$. By definition, $\underline\delta(p)=1/(dp)$. Lemma~\ref{lem:giveup-psi} gives $\bar\delta(p)\ge1/(dq)$, and $q\le p$ therefore implies $\underline\delta(p)\le\bar\delta(p)$. Suppose $\bar\delta(p)<\infty$ and $\delta\ge\bar\delta(p)$. By definition, $\bar\delta(p)$ is the smallest value at which $\Psi_q\le 0$. By monotonicity, $\Psi_q(\delta)\le\Psi_q(\bar\delta(p))\le0$.
Thus Lemma~\ref{lem:giveup-lower} makes the zero estimator minimax-optimal for $\delta\ge\bar\delta(p)$. If $\delta < \bar\delta(p)$, $\Psi_q(\delta)>0$ and Lemma~\ref{lem:giveup-lower} rules out zero as a minimax-optimal estimator.
Applying Lemma~\ref{lem:giveup-psi} with $r=q$ also gives the finiteness condition $\bar\delta(p)<\infty$ if and only if $p^d\le1/2$.

\emph{Excluding control outcomes.} Fix $p>1/2$ and set $C_p=12M_q$. Choose $N_p$ large enough that, whenever $n\ge N_p$,
\[
c\coloneqq\frac{1+C_p/n}{p}<2,\qquad
n>\frac{8qM_q^2}{p}.
\]
Fix such an $n$ and $\delta\ge c/d$. The first part of the proof already covers $\delta\ge\bar\delta(p)$, so assume $\Psi_q(\delta)>0$. Because outcomes lie in $[-1,1]$, we may replace $\delta$ by $\min\{\delta,2\}$ without changing the outcome class. Since $c/d<2$, the condition $\delta\ge c/d$ is preserved. We may therefore assume $\delta\le2$.

Let $a$ minimize $F(0,t)$ over $t\in\mathbb R$. We will show that $(0,a)$ remains optimal when the control coefficient can also vary. Since $a$ minimizes $F(0,t)$, the directional derivatives of $F$ at $(0,a)$ in the directions $(0,1)$ and $(0,-1)$ are nonnegative. By~\eqref{eq:giveup-directional}, there are two outcome functions $y^{\uparrow},y^{\downarrow}\in\mathcal Y^\star(0,a)$ whose treatment derivatives satisfy
\[
\partial_{a^1}\mathcal R\big((0,a),y^{\uparrow}\big)\ge0\ge\partial_{a^1}\mathcal R\big((0,a),y^{\downarrow}\big).
\]
With treatment outcomes fixed and $a^0=0$, risk is a strictly convex quadratic in $\sum_i y_i(\mathbf0)\in[-n,n]$. Thus each worst-case outcome function has either $y_i(\mathbf0)=1$ for every $i$ or $y_i(\mathbf0)=-1$ for every $i$. Reversing all outcome signs preserves the risk at every coefficient pair, and hence its derivatives, so take both $y^\uparrow$ and $y^\downarrow$ to satisfy $y_i(\mathbf0)=1$ for every $i$.

Now fix a direction $\mathbf w=(w^0,w^1)$. First select the base outcome function according to the treatment component:
\[
y\coloneqq
\begin{cases}
y^{\uparrow},&w^1\ge0,\\
y^{\downarrow},&w^1<0.
\end{cases}
\]
By construction, $w^1\,\partial_{a^1}\mathcal R((0,a),y)\ge0$. Next apply the constant or linear control construction of Lemma~\ref{lem:giveup-control-profiles} to the selected $y$, according to the control component. Explicitly, for every $i$ and $\mathbf z$, set
\[
y'_i(\mathbf z)\coloneqq
\begin{cases}
z_i y_i(\mathbf z)+(1-z_i),&w^0\ge0,\\
z_i y_i(\mathbf z)+(1-z_i)\left(1-\dfrac{c\,s_i(\mathbf z)}d\right),&w^0<0.
\end{cases}
\]
Thus $y'$ retains the treatment outcomes of $y^{\uparrow}$ when $w^1\ge0$ and of $y^{\downarrow}$ when $w^1<0$; the sign of $w^0$ determines its control outcomes. By Lemma~\ref{lem:giveup-control-profiles} and the equality of treatment derivatives established in its proof, $y'\in\mathcal Y^\star(0,a)$ and
\[
\begin{aligned}
w^1\,\partial_{a^1}\mathcal R\big((0,a),y'\big)
&=w^1\,\partial_{a^1}\mathcal R\big((0,a),y\big)\ge0,\\
w^0\,\partial_{a^0}\mathcal R\big((0,a),y'\big)&\ge0.
\end{aligned}
\]
Hence
\[
\nabla\mathcal R\big((0,a),y'\big)^{\mathsf T}\mathbf w
=w^0\,\partial_{a^0}\mathcal R\big((0,a),y'\big)+w^1\,\partial_{a^1}\mathcal R\big((0,a),y'\big)\ge0.
\]
Since $\mathbf w$ was arbitrary, the optimality condition in Lemma~\ref{lem:giveup-danskin} shows that $(0,a)$ minimizes $F$, i.e., $(0,a)$ is an optimal uniform estimator. The symmetry reduction in Theorem~\ref{thm:sym-static-general-p} makes it optimal among all static estimators as well. The argument holds for every $\delta\ge c/d$, and exchanging treatment and control proves the case $p<1/2$.
\end{proof}

\subsection{Optimal risk across interference levels}\label{app:risk-curves}

\begin{proof}[Proof of Theorem~\ref{thm:risk-bounds}]
Write $q=1-p$. By Theorem~\ref{thm:sym-static-general-p}, we may use the same coefficient $a^1$ for every treated unit and the same coefficient $a^0$ for every control unit. We first choose an estimator and bound its bias and variance, then prove the lower bound for arbitrary coefficients.

\textbf{Proof of the upper bound.}
Choose a uniform estimator $a=(a^0,a^1)$ with
\[
a^1=\frac1p\mathbf1\{dq\delta<1\},\qquad
a^0=-\frac1q\mathbf1\{dp\delta<1\}.
\]
For any allowed outcome function $y$, its bias is
\begin{equation}\label{eq:risk-estimator-bias}
\begin{aligned}
\mathbb E[\mathbf a(\mathbf Z)^{\mathsf T}\mathbf y(\mathbf Z)]-GTE(y)
&=\sum_{i=1}^n\left(pa^1\mathbb E[y_i(\mathbf Z)\mid Z_i=1]-y_i(\mathbf1)\right)+\sum_{i=1}^n\left(qa^0\mathbb E[y_i(\mathbf Z)\mid Z_i=0]+y_i(\mathbf0)\right).
\end{aligned}
\end{equation}
Bounded interference limits the change in unit $i$'s outcome by $\delta$ times the number of neighbors whose assignments change. Independence of the assignments therefore gives
\[
\begin{aligned}
\left|\mathbb E[y_i(\mathbf Z)\mid Z_i=1]-y_i(\mathbf1)\right|
&\le\delta\,\mathbb E[d-s_i(\mathbf Z)\mid Z_i=1]=dq\delta,\\
\left|\mathbb E[y_i(\mathbf Z)\mid Z_i=0]-y_i(\mathbf0)\right|
&\le\delta\,\mathbb E[s_i(\mathbf Z)\mid Z_i=0]=dp\delta.
\end{aligned}
\]
Add and subtract $pa^1y_i(\mathbf1)$ in the first sum of~\eqref{eq:risk-estimator-bias}, and $qa^0y_i(\mathbf0)$ in the second. Using $|y_i(\mathbf1)|,|y_i(\mathbf0)|\le1$ and the triangle inequality yields
\[
\begin{aligned}
|\mathbb E[\mathbf a(\mathbf Z)^{\mathsf T}\mathbf y(\mathbf Z)]-GTE(y)|
&\le n\left(|pa^1|dq\delta+|pa^1-1|+|qa^0|dp\delta+|qa^0+1|\right)\\
&=n\big[\min\{dq\delta,1\}+\min\{dp\delta,1\}\big]
=n\sqrt{\Phi_{d,p}(\delta)}.
\end{aligned}
\]
Here the equality follows by substituting the chosen coefficients: $pa^1$ is 1 when $dq\delta<1$ and 0 otherwise, while $qa^0$ is -1 when $dp\delta<1$ and 0 otherwise.

We next bound the variance. Changing one unit's assignment changes its contribution to the estimator by at most $|a^0|+|a^1|$. Only its $d$ neighbors' outcomes can also change; each changes by at most $\delta$, and its coefficient has magnitude at most $\max\{|a^0|,|a^1|\}$. The chosen coefficients satisfy
\[
\begin{aligned}
d\delta|a^1|&=\frac{dq\delta\,\mathbf1\{dq\delta<1\}}{pq}\le\frac1{pq},\\
d\delta|a^0|&=\frac{dp\delta\,\mathbf1\{dp\delta<1\}}{pq}\le\frac1{pq}.
\end{aligned}
\]
Consequently, changing one assignment changes the estimator by at most
\[
|a^0|+|a^1|+d\delta\max\{|a^0|,|a^1|\}
\le\frac1q+\frac1p+\frac1{pq}=\frac2{pq}.
\]
Lemma~\ref{lem:giveup-variance} therefore gives
\[
\operatorname{Var}(\mathbf a(\mathbf Z)^{\mathsf T}\mathbf y(\mathbf Z))\le npq\left(\frac2{pq}\right)^2=\frac{4n}{pq}.
\]
Combining the bias and variance bounds, for every allowed $y$,
\[
\mathbb E(\mathbf a(\mathbf Z)^{\mathsf T}\mathbf y(\mathbf Z)-GTE(y))^2=(\mathbb E[\mathbf a(\mathbf Z)^{\mathsf T}\mathbf y(\mathbf Z)]-GTE(y))^2+\operatorname{Var}(\mathbf a(\mathbf Z)^{\mathsf T}\mathbf y(\mathbf Z))
\le n^2\Phi_{d,p}(\delta)+\frac{4n}{pq}.
\]
Taking the worst case over $y$ and dividing by $n^2$ proves the upper bound with error $4/(npq)$. The zero estimator also gives risk $4n^2$, completing the upper bound in~\eqref{eq:risk-bounds}.

\textbf{Proof of the lower bound.}
Fix $(a^0,a^1)$. Mean squared error is at least squared bias.

For $r\in\{p,q\}$, set
\[
\rho_r=\min\{\delta,1/(dr)\},\qquad
c_r=\max\left\{0,dr\rho_r-\frac{dr(1-r)\rho_r^2}{4}\right\}.
\]
For every unit $i$ and assignment $\mathbf z$, define the outcome function by
\[
y_i^{a^0,a^1}(\mathbf z)=
\begin{cases}
-1, & z_i=1,\ a^1<0,\\
\min\{1,-c_q+\rho_q[d-s_i(\mathbf z)]\}, & z_i=1,\ a^1\ge0,\\
1, & z_i=0,\ a^0>0,\\
-\min\{1,-c_p+\rho_p s_i(\mathbf z)\}, & z_i=0,\ a^0\le0.
\end{cases}
\]
The constant cases lie in $[-1,1]$. In the other cases, $0\le c_r\le dr\rho_r\le1$, so $-c_r+\rho_r k\ge-1$ for $0\le k\le d$; clipping at 1 and changing the sign preserve the bound. For $j\ne i$, changing $z_j$ leaves $z_i$ fixed. If $j$ is not a neighbor, $s_i(\mathbf z)$ does not change. If $j$ is a neighbor, $s_i(\mathbf z)$ changes by one, and $x\mapsto\min\{1,x\}$ is 1-Lipschitz. The change in $y_i^{a^0,a^1}$ is therefore zero in a constant case, at most $\rho_q$ in the treated clipped case, and at most $\rho_p$ in the control clipped case. Hence
\[
\left|y_i^{a^0,a^1}(\mathbf z^{(z_j=1)})-y_i^{a^0,a^1}(\mathbf z^{(z_j=0)})\right|
\le\max\{\rho_p,\rho_q\}\mathbf A_{ij}\le\delta\mathbf A_{ij}.
\]
Thus $y^{a^0,a^1}$ satisfies both outcome constraints. No constraint links the cases $z_i=1$ and $z_i=0$, because bounded interference only concerns $j\ne i$.

Let $B_r\sim\operatorname{Bin}(d,r)$. Given $Z_i=1$, the number of control neighbors has distribution $B_q$; given $Z_i=0$, the number of treated neighbors has distribution $B_p$. Substitution into~\eqref{eq:risk-estimator-bias} gives the two terms for each unit:
\[
\begin{aligned}
&pa^1\mathbb E[y_i(\mathbf Z)\mid Z_i=1]-y_i(\mathbf1)\\
&\qquad=\begin{cases}
1-pa^1,&a^1<0,\\
c_q+pa^1\mathbb E\min\{1,-c_q+\rho_qB_q\},&a^1\ge0,
\end{cases}\\[3pt]
&qa^0\mathbb E[y_i(\mathbf Z)\mid Z_i=0]+y_i(\mathbf0)\\
&\qquad=\begin{cases}
1+qa^0,&a^0>0,\\
c_p-qa^0\mathbb E\min\{1,-c_p+\rho_pB_p\},&a^0\le0.
\end{cases}
\end{aligned}
\]

The expectations in the second cases are nonnegative. For $r\in\{p,q\}$, the claim is immediate when $c_r=0$. When $c_r>0$, put $m=dr\rho_r-c_r=dr(1-r)\rho_r^2/4$. The identity $\min\{1,z\}=(1+z-|z-1|)/2$, the variance $\operatorname{Var}(\rho_rB_r)=4m$, and Cauchy--Schwarz give
\[
\mathbb E\min\{1,-c_r+\rho_rB_r\}
\ge\frac{1+m-\sqrt{(m-1)^2+4m}}2=0.
\]
For $a^1<0$, the first treatment case is $1-pa^1\ge1\ge c_q$; for $a^0>0$, the first control case is $1+qa^0\ge1\ge c_p$. In the remaining cases, $pa^1\ge0$ and $-qa^0\ge0$, so the nonnegative expectations give a treatment term of at least $c_q$ and a control term of at least $c_p$. Summing over units gives $\mathbb E[\mathbf a(\mathbf Z)^{\mathsf T}\mathbf y(\mathbf Z)]-GTE(y)\ge n(c_p+c_q)$. Taking the minimum over $(a^0,a^1)$ and using the symmetry reduction at the start of the proof gives
\[
\frac{V_{G,p}(\delta)}{n^2}\ge(c_p+c_q)^2.
\]
Finally, $dr\rho_r=\min\{dr\delta,1\}$ and $\rho_r\le1/(dr)$ imply
\[
0\le\min\{dr\delta,1\}-c_r\le\frac{1-r}{4dr}.
\]
Both $c_p+c_q$ and $\sqrt{\Phi_{d,p}(\delta)}$ lie in $[0,2]$. Taking the difference of their squares therefore gives
\[
\Phi_{d,p}(\delta)-(c_p+c_q)^2
\le4\left(\frac{q}{4dp}+\frac{p}{4dq}\right)
=\frac{p^2+q^2}{dpq}.
\]
Combining the last inequality with nonnegativity of risk proves the lower bound.
\end{proof}

%% file: references.bib
@book{bertsekas1999nonlinear,
  author    = {Bertsekas, Dimitri P.},
  title     = {Nonlinear Programming},
  edition   = {2nd},
  publisher = {Athena Scientific},
  address   = {Belmont, MA},
  year      = {1999}
}

@book{boucheron2013concentration,
  author    = {Boucheron, St{\'e}phane and Lugosi, G{\'a}bor and Massart, Pascal},
  title     = {Concentration Inequalities: A Nonasymptotic Theory of Independence},
  publisher = {Oxford University Press},
  address   = {Oxford},
  year      = {2013}
}

@article{holtz2025reducing,
  title={Reducing interference bias in online marketplace experiments using cluster randomization: Evidence from a pricing meta-experiment on airbnb},
  author={Holtz, David and Lobel, Felipe and Lobel, Ruben and Liskovich, Inessa and Aral, Sinan},
  journal={Management Science},
  volume={71},
  number={1},
  pages={390--406},
  year={2025},
  publisher={INFORMS}
}

@article{bojinov2023design,
  title={Design and analysis of switchback experiments},
  author={Bojinov, Iavor and Simchi-Levi, David and Zhao, Jinglong},
  journal={Management Science},
  volume={69},
  number={7},
  pages={3759--3777},
  year={2023},
  publisher={INFORMS}
}

@article{ni2023design,
  title={Design of panel experiments with spatial and temporal interference},
  author={Ni, Tu and Bojinov, Iavor and Zhao, Jinglong},
  journal={Available at SSRN 4466598},
  year={2023}
}

@article{boutilier2024operational,
  title={Operational dosage: Implications of capacity constraints for the design and interpretation of experiments},
  author={Boutilier, Justin and Jonasson, Jonas Oddur and Li, Hannah and Yoeli, Erez},
  journal={arXiv preprint arXiv:2407.21322},
  year={2024}
}

@article{johari2024does,
  title={When does interference matter? decision-making in platform experiments},
  author={Johari, Ramesh and Li, Hannah and Murthy, Anushka and Weintraub, Gabriel Y},
  journal={arXiv preprint arXiv:2410.06580},
  year={2024}
}

@article{han2024population,
  title={Population interference in panel experiments},
  author={Han, Kevin and Basse, Guillaume and Bojinov, Iavor},
  journal={Journal of Econometrics},
  volume={238},
  number={1},
  pages={105565},
  year={2024},
  publisher={Elsevier}
}

@article{peng2025differences,
  title={Differences-in-neighbors for network interference in experiments},
  author={Peng, Tianyi and Ye, Naimeng and Zheng, Andrew},
  journal={arXiv preprint arXiv:2503.02271},
  year={2025}
}

@article{farias2022markovian,
  title={Markovian interference in experiments},
  author={Farias, Vivek and Li, Andrew and Peng, Tianyi and Zheng, Andrew},
  journal={Advances in Neural Information Processing Systems},
  volume={35},
  pages={535--549},
  year={2022}
}

@book{ben2009robust,
  title={Robust optimization},
  author={Ben-Tal, Aharon and El Ghaoui, Laurent and Nemirovski, Arkadi},
  year={2009},
  publisher={Princeton university press}
}

@book{bertsimas2022robust,
  title={Robust and adaptive optimization},
  author={Bertsimas, Dimitris and {den Hertog}, Dick},
  publisher={Dynamic Ideas},
  year={2022}
}

@article{blanchet2025distributionally,
  title={Distributionally robust optimization and robust statistics},
  author={Blanchet, Jose and Li, Jiajin and Lin, Sirui and Zhang, Xuhui},
  journal={Statistical Science},
  volume={40},
  number={3},
  pages={351--377},
  year={2025},
  publisher={Institute of Mathematical Statistics}
}

@article{bright2025reducing,
  title={Reducing marketplace interference bias via shadow prices},
  author={Bright, Ido and Delarue, Arthur and Lobel, Ilan},
  journal={Management Science},
  volume={71},
  number={8},
  pages={7094--7112},
  year={2025},
  publisher={INFORMS}
}

@article{chen2025b,
  title={A/b testing under product stockouts},
  author={Chen, Liting and Cohen, Maxime C and Miao, Sentao},
  journal={Available at SSRN 5256420},
  year={2025}
}

@article{li2026experimenting,
  title={Experimenting under stochastic congestion},
  author={Li, Shuangning and Johari, Ramesh and Kuang, Xu and Wager, Stefan},
  journal={Management Science},
  year={2026},
  publisher={INFORMS}
}

@inproceedings{ugander2013graph,
  title={Graph cluster randomization: Network exposure to multiple universes},
  author={Ugander, Johan and Karrer, Brian and Backstrom, Lars and Kleinberg, Jon},
  booktitle={Proceedings of the 19th ACM SIGKDD international conference on Knowledge discovery and data mining},
  pages={329--337},
  year={2013}
}

@article{jia2023clustered,
  title={Clustered switchback designs for experimentation under spatio-temporal interference},
  author={Jia, Su and Kallus, Nathan and Yu, Christina Lee},
  journal={arXiv preprint arXiv:2312.15574},
  year={2023}
}

@article{aronow2017estimating,
  title={Estimating average causal effects under general interference, with application to a social network experiment},
  author={Aronow, Peter M and Samii, Cyrus},
  year={2017}
}

@article{brennan2022cluster,
  title={Cluster randomized designs for one-sided bipartite experiments},
  author={Brennan, Jennifer and Mirrokni, Vahab and Pouget-Abadie, Jean},
  journal={Advances in Neural Information Processing Systems},
  volume={35},
  pages={37962--37974},
  year={2022}
}

@article{candogan2024correlated,
  title={Correlated cluster-based randomized experiments: Robust variance minimization},
  author={Candogan, Ozan and Chen, Chen and Niazadeh, Rad},
  journal={Management Science},
  volume={70},
  number={6},
  pages={4069--4086},
  year={2024},
  publisher={INFORMS}
}

@article{viviano2023causal,
  title={Causal clustering: design of cluster experiments under network interference},
  author={Viviano, Davide and Lei, Lihua and Imbens, Guido and Karrer, Brian and Schrijvers, Okke and Shi, Liang},
  journal={arXiv preprint arXiv:2310.14983},
  year={2023}
}

@article{kandiros2026design,
  title={A Design-Based Minimax Theory for Network Experiments},
  author={Kandiros, Vardis and Harshaw, Christopher and S{\"a}vje, Fredrik},
  journal={arXiv preprint arXiv:2608.04909},
  year={2026}
}

@article{thiyageswaran2026optimal,
  title={Optimal Design under Interference, Homophily, and Robustness Trade-offs},
  author={Thiyageswaran, Vydhourie and Kokot, Alex and Brennan, Jennifer and Meila, Marina and Yu, Christina Lee and Fazel, Maryam},
  journal={arXiv preprint arXiv:2601.17145},
  year={2026}
}

@article{sudijono2026sharp,
  title={Sharp Minimax Theory for Randomized Experiments},
  author={Sudijono, Timothy and Dobriban, Edgar and Tchetgen, Eric Tchetgen},
  journal={arXiv preprint arXiv:2608.13822},
  year={2026}
}
